\documentstyle {article}
\newcommand{\vsp}{\medskip} 
\newcommand{\be}{\begin{equation}}
\newcommand{\ee}{\end{equation}}
\newcommand{\bea}{\begin{eqnarray}}
\newcommand{\eea}{\end{eqnarray}}
\newcommand{\nn}{\nonumber}

\newcommand{\tC}{\hat{C}(A)=\oplus_{i,j\geq 0}Hom(A^{\ot i};A^{\ot j})}
\newcommand{\bC}{\bar{C}(A)=\oplus_{j\geq 0}Hom(A;A^{\ot j})}
\newcommand{\bde}{\bar{\delta}}
\newcommand{\hde}{\hat{\delta}}
\newcommand{\ma}{{\rm max}}
\newcommand{\com}{\Delta}
\newcommand{\tx}{\tilde{x}}
\newcommand{\mxy}{\{ m\}\{\{ x_1\}\{ y_1,\dots,y_k\},\{ x_2\}\{ y_{k+1},
    \dots,y_n\}\}}
\newcommand{\m}{m}
\newcommand{\sa}{sa}
\newcommand{\CA}{C^{\bullet}(A)}
\newcommand{\CB}{C^{\bullet}(B)}

\newcommand{\CBB}{C^{\bullet,\bullet}(A)}
\newcommand{\HTT}{Hom(TA;TA)}

\newcommand{\tm}{\tilde{m}}

\newcommand{\lbr}{\lceil}
\newcommand{\rbr}{\rceil}
\newcommand{\de}{\delta}
\newcommand{\De}{\Delta}
\newcommand{\Ai}{A_{\infty}}
\newcommand{\Li}{L_{\infty}}
\newcommand{\op}{\oplus}
\newcommand{\ot}{\otimes}
\newcommand{\ra}{\rightarrow}
\newcommand{\Z}{{\bf Z}}
\newcommand{\C}{{\bf C}}
\newcommand{\oper}[3]{\Phi^{#2}\mbox{\bf [}#3;#1\mbox{\bf ]}}
\newcommand{\defi}{{\stackrel{\rm def}{=}}}
\newcommand{\po}{\sum_{t=1}^{n}(n-t)\| a_t\|}
\newcommand{\pos}{\sum_{t=1}^{n}(n-t)\| a_{\sigma(t)}\|}
\newcommand{\pot}{\sum_{t=1}^{n}(n-t)\| a_{\tau(t)}\|}

\newcommand{\p}{p(\sigma ;a_1,\dots,a_n)}
\newcommand{\pttt}{\tilde{\epsilon}(\sigma ;a_1,\dots,a_n)}
\newcommand{\su}{\sum_{u<v,\sigma^{-1}(u)>\sigma^{-1}(v)}
(\| a_u\| +\| a_v\| )}
\newcommand{\sut}{\sum_{u<v,\tau^{-1}(u)>\tau^{-1}(v)}
(\| a_u\| +\| a_v\| )}
\newcommand{\adj}{{\rm ad}}
\newcommand{\ide}{{\rm id}}
\newcommand{\ep}{\epsilon}
\newcommand{\tep}{\tilde{e}}
\newcommand{\uv}{u<v,\sigma^{-1}(u)>\sigma^{-1}(v)}
\newcommand{\dep}{\delta'}
\newcommand{\bj}{\beta_j}
\newcommand{\gi}{\gamma_{i-1}}
\newcommand{\gak}{\gamma_k}
\newcommand{\gakb}{\gamma_{k+1}}
\newcommand{\dbb}{\de(\beta_1)}
\newcommand{\dbj}{\de(\bj)}
\newcommand{\dgb}{\dep(\gamma_1)}
\newcommand{\dgi}{\dep(\gi)}
\newcommand{\dgk}{\dep(\gak)}
\newcommand{\dgkb}{\dep(\gakb)}
\newcommand{\si}{\sigma}
\newcommand{\sj}{\si_J}
\newcommand{\tes}{\tep(\si ;a_1,\dots,a_n)}
\newcommand{\tessu}{\tep(\sigma_3 ;a_1,\dots,a_n)}
\newcommand{\tesd}{\tep(\de ;a_{\beta_1},\dots,a_{\bj})}
\newcommand{\tesdp}{\tep(\dep ;a_{\gamma_1},\dots,a_{\gi})}
\newcommand{\susb}{\sum_{t=1}^k\| a_{\dep(\gamma_t)}\|\sum_{T=1}^j\| 
    a_{\beta_T}\|}
\newcommand{\susa}{\sum_{t=1}^k\| a_{\dep(\gamma_t)}\|\sum_{T=1}^j\|
    a_{\de(\beta_T)}\|} 
\newcommand{\tebiu}{\tep(\sigma_1\sigma_2\sigma_3 ;a_1,\dots,a_n)}
\newcommand{\tessi}{\tep(\sigma_2 ;a_{\beta_1},\dots,a_{\bj},a_{\gamma_1},\dots,
    a_{\gi})}
\newcommand{\tessb}{\tep(\sigma_1 ;a_{\de(\beta_1)},\dots,a_{\de(\bj)},
    a_{\dep(\gamma_1)},\dots,a_{\dep(\gi)})}
\newcommand{\dpjp}{\dep\in{\rm Perm}(J')}
\newcommand{\dpj}{\de\in{\rm Perm}(J)}
\newcommand{\tessd}{\tep(\sigma_4 ;a_0,a_{\gamma_1},\dots,a_{\gi})}
\newcommand{\susi}{\| a_0\|\sum_{t=1}^k\| a_{\dep(\gamma_t)}\|}
\newcommand{\mms}{\{\tilde{m}\circ\tilde{m}\}\{ sa_1\}\cdots\{ sa_n\}}
\newcommand{\tej}{\tilde{\epsilon}(\sj ;a_1,\dots,a_n)}
\newcommand{\lilj}{\tilde{l}_i(\tilde{l}_j(a_{\beta_1},\dots,a_{\bj}),
    a_{\gamma_1},\dots,a_{\gi})}
\newcommand{\susu}{\sum_{t=1}^k\| a_{\si(t)}\|}
\newcommand{\mimj}{\tilde{m}_i(a_{\si(1)},\dots,a_{\si(k)},\tilde{m}_j
     (a_{\si(k+1)},\dots,a_{\si(k+j)}),a_{\si(k+j+1)},\dots,a_{\si(n)})} 
\newcommand{\susd}{\sum_{t=1}^k\| a_{\dep(\gamma_t)}\|(\sum_{T=1}^j\|
     a_{\beta_T}\| +1)}
\newcommand{\mimjt}{\tilde{m}_i(a_{\dep(\gamma_1)},\dots,a_{\dep(\gak)},
     \tilde{m}_j(a_{\de(\beta_1)},\dots,a_{\de(\bj)}),a_{\dep(\gakb)},\dots,
     a_{\dep(\gi)})}
\newcommand{\milj}{\tilde{m}_i(a_{\dep(\gamma_1)},\dots,a_{\dep(\gak)},
     \tilde{l}_j(a_{\beta_1},\dots,a_{\bj}),a_{\dep(\gakb)},\dots,
     a_{\dep(\gi)})}
\newcommand{\gki}{\gamma_{i-1}}
\newcommand{\mD}{{[m,\tri ]}}
\newcommand{\lls}{\{\tilde{l}\circ\tilde{l}\}\{ sa_1\}\cdots\{ sa_n\}} 
\newcommand{\LiLj}{\tilde{l}_i(\tilde{l}_j(a_{\beta_1},\dots,a_{\beta_j}),
      a_{\gamma_1},\dots,a_{\gamma_{i-1}})}
\newcommand{\LiLjt}{l_i(l_j(a_{\beta_1},\dots,a_{\beta_j}),a_{\gamma_1},
      \dots,a_{\gamma_{i-1}})}
\newcommand{\liljt}{\tilde{l}_i(a_{\dep(\gamma_1)},\dots,a_{\dep(\gamma_k)},
      \tilde{l}_j(a_{\de(\beta_1)},\dots,a_{\de(\beta_j)}),
      a_{\dep(\gamma_{k+1})},\dots,a_{\dep(\gamma_{i-1})})}
\newcommand{\lilju}{\tilde{l}_i(a_{\dep(\gamma_1)},\dots,a_{\dep(\gamma_k)},
      \tilde{l}_j(a_{\beta_1},\dots,a_{\beta_j}),a_{\dep(\gamma_{k+1})},
      \dots,a_{\dep(\gamma_{i-1})})}
\newcommand{\sumj}{\sum_{u<v,\sigma_J^{-1}(u)>\sigma_J^{-1}(v)}\| a_u\|\,
      \| a_v\|} 
\newcommand{\sumjbt}{\sum_{t=1}^j(j-t)\| a_{\beta_t}\|}
\newcommand{\sumbt}{\sum_{t=1}^j\| a_{\beta_t}\|}
\newcommand{\sumig}{\sum_{t=1}^{i-1}(i-1-t)\| a_{\gamma_t}\|}
\newcommand{\tri}{\triangle}

\newtheorem{lemma}{Lemma}
\newtheorem{thm}{Theorem}
\newtheorem{prop}{Proposition}
\newtheorem{ex}{Example}
\newtheorem{Rm}{Remark}
\newtheorem{defn}{Definition}
\newtheorem{cor}{Corollary}

\begin{document}

\title{{\bf Multibraces on the Hochschild complex}}
\author{F\"{U}SUN AKMAN\thanks{ Partially presented at the AMS meeting in
Lawrenceville, NJ, Oct. 5-6, 1996. Part of this work was done at Cornell
University.} \\
Dept. of Mathematics and Statistics\\ Utah State University \\ Logan, UT
84322-3900 \\ fusun@math.usu.edu}
\date{August 15, 1997}
\maketitle

\begin{abstract}
We generalize the coupled braces $\{ x\}\{ y\}$ of Gerstenhaber and $\{ x\}
\{ y_1,\dots,y_n\}$ of Getzler depicting compositions of multilinear
maps in the Hochschild complex $\CA =Hom(TA;A)$
of a graded vector space $A$ to expressions of the form $\{ x_1^{(1)},\dots,
 x_{i_1}^{(1)}\}\cdots\{ x_{1}^{(m)},\dots,x_{i_m}^{(m)}\}$ on the extended
space $\CBB =\HTT$, and clarify many of the existing sign conventions that
show up in the algebra of mathematical physics (namely in associative and
Lie algebras, Batalin-Vilkovisky algebras, $\Ai$ and $\Li$ algebras).
As a result, we introduce a new variant of the master identity for $\Li$
algebras. We also 
comment on the bialgebra cohomology differential of Gerstenhaber and Schack,
and define 
multilinear higher order differential operators with respect to multilinear
maps using the new language. The continuation of this work will be on the
various homotopy structures on a topological vertex operator algebra,
as introduced by Kimura, Voronov, and Zuckerman.
\end{abstract}

\tableofcontents

\section {Introduction}

Mathematical physics and homological algebra are infested with magical
algebraic identities which usually boil down to the following: the
composition of a multilinear map with another, or a sum of such
compositions, is identically zero. The lack of a unifying language makes it
hard to see the origins and generalizations of such statements, as well as
to prove them. Many examples will be shown among explicit formulas and
properties of differentials in cohomology theories, higher homotopy
algebras, and especially algebraic identities arising from topological
operads in mathematical physics; the author's interest in the subject began
with Kimura, Voronov, and Zuckerman's ``Homotopy Gerstenhaber algebras and
topological field theory'' \cite{KVZ}, where multilinear expressions
\[ \{ v_1,\dots,v_m\}\cdots\{ w_1,\cdots,w_n\} \]
(arguments living in a topological vertex operator algebra, or TVOA)
satisfy some lower identities resembling those for the braces $\{ x\}\{
y\}$ of Gerstenhaber (from the 1960's!) and $\{ x\}\{ y_1,\dots, y_n\}$ of
Getzler, which denote the substitution of the multilinear map(s) on the right
into the one on the left. Unlike the braces in~\cite{KVZ}, the latter two
did not extend beyond two pairs, except in iterations, and it seemed
natural to stretch the idea as far as possible since the literature was now
ripe for new usage. It turns out that multibraces are indeed both a
convenient language and a shortcut for {\sl expressing} many ideas, their
usefulness readily demonstrated in {\sl proving} and {\sl generalizing}
statements. In particular, we will introduce a new master identity for
strongly homotopy Lie algebras which are obtained by antisymmetrizing
products in strongly homotopy associative algebras (Theorem~\ref{equi}), 
and point
out more general definitions of many concepts, such as higher order
differential operators on noncommutative, nonassociative algebras 
(Section~\ref{phi}). Simple proofs of new and old results will make heavy
use of the multibraces language. Here is an outline of this paper:
\vsp

The {\em coupled pairs of braces} 
\[ \{ x\} \{ y\} =x\circ y \]
of Gerstenhaber, and
\be\label{get} \{ x\} \{ y_1,\dots,y_n\} \ee
of Getzler on the Hochschild complex
\[ \CA =Hom(TA;A) \]
of a graded vector space $A$, were defined to be generalizations of 
substitution of
elements of $A$ into a multilinear map, and of composition of linear 
maps on $A$, where order and grading are extremely important (the first
pair of braces from the left were omitted in the original works; we adopt
the uniform notation of Kimura, Voronov and Zuckerman in \cite{KVZ}). We 
will go over their definitions and propose yet another generalization
\be \{ x\} \{ x_1^{(1)},\dots,x_{i_1}^{(1)}\} \cdots \{ x_1^{(m)},\dots,
x_{i_m}^{(m)} \} \{ a_1,\dots,a_n\} \label{exp} \ee
of this formalism, where incomplete expressions
\be \{ x\} \{ x_1^{(1)},\dots,x_{i_1}^{(1)} \} \cdots \{ x_1^{(m)},\dots,
x_{i_m}^{(m)} \}\label{ts}\ee
(i.e. those which haven't been fed some $a_1$, ..., $a_n$)
are understood to be multilinear maps with values in $A$, which are 
eventually evaluated at
\[ \{ a_1,\dots,a_n\} ,\]
or even at
\[ \{ a_1,\dots\} \cdots \{ \dots,a_n\} , \;\; a_i\in A. \]
In its simplest form, 
\[ \{ x\}\{ y\}\{ a\} =x(y(a)) \]
is the substitution of $a$ into the composition $x\circ y$ of linear functions
$x$ and $y$ on $A$.
Such expressions preserve the {\em (adjusted) degree of homogeneity} $d$ of
elements of
$\CA$ ($d(x)=n-1$ if $x$ is $n$-linear; $d(a)=-1$ if $a\in A$). The definition
of (\ref{exp}) is of the ``follow your nose'' variety, as a result of which
an expression like 
\[ \{ m\} \{ a_1\} \cdots \{ a_n\} \;\;\;\;\;\; (d(m)=n-1,\; d(a_i)=-1) \]
stands for a (signed) sum over permutations of
\[ \{ m\} \{ a_1,\dots,a_n\} \stackrel{\rm def}{=} m(a_1,\dots,a_n).\]
In general,
\[ \{ m\} \{ a_1,\dots\} \cdots\{ \dots,a_n\} \]
is just a sum over those permutations which fix the order within each 
individual string of $a_i$'s. Then two expressions of type~(\ref{ts}) are 
deemed equal as multilinear maps if they are equal when evaluated at 
all $\{ a_1,\dots,a_n\}$. We emphasize that although~(\ref{exp}) can be
obtained by multiple iterations of~(\ref{get}), the full potential of
the formalism in~(\ref{get}), especially regarding substitution, has not
been attained
so far. Moreover, coupled pairs of braces can be used to denote elements
of the tensor algebra $TA$ and operations with values in $TA$. As expected,
$\{ a_1,\dots,a_n\}$ stands for $a_1\ot\cdots\ot a_n\in A^{\ot n}$, 
$\{\{ a_1,\dots,a_k\} ,\{ a_{k+1},\dots,a_n\}\}'$ stands for $(a_1\ot\cdots
\ot a_k)\ot(a_{k+1}\ot\cdots\ot a_n)\in A^{\ot k}\ot A^{\ot(n-k)}$ (as opposed
to $A^{\ot n}$; primed braces belong to $T(TA)$ by definition),
and the symbol $\{ x,y\}\{ a_1,\dots,a_n\}$ will mean $\pm
x(a_1,\dots,a_k)\ot y(a_{k+1},\dots,a_n)\in TA$ for appropriate $x$ and $y$.
In short, we will expand our notation to
\[\{ x_1^{(1)},\dots,x_{i_1}^{(1)}\}\cdots\{ x_1^{(m)},\dots,x_{i_m}^{(m)}\}\]
on
\[ \CBB =\HTT\]
in Section~\ref{twofour} (until then, all our maps will have values in A).
In particular, the coupled-braces notation can be used to write explicit
practical formulas for the bialgebra cohomology differential of Gerstenhaber 
and Schack \cite{GeS}.
\vsp

This language makes many concepts and proofs very accessible in multilinear 
algebra
(see the proof of $\de^2=0$ on the complex $\CA$ for an associative 
or $\Ai$ algebra and of various BV algebra identities in Sections \ref{AA},
\ref{BB}, and \ref{CC} respectively). 
For example, an associative algebra is just a graded vector space $A$ 
with $m\in C^2(A)$ satisfying
\[ m\circ m=0; \]
one can make a Lie algebra out of it via the brackets
\[ \lbr a,b\rbr =\{ m\} \{ a\} \{ b\} .\]
An $\Ai$ algebra is again some $A$ with $m\in\CA$,
\[ m=m_1+m_2+\cdots ,\;\;\; d(m_k)=k-1,\;\;\; (-1)^{|m_k|}=(-1)^k,\]
satisfying the master identity
\[ \tilde{m}\circ\tilde{m}=0,\;\;\;\mbox{or}\;\;\;\{\tilde{m}\circ\tilde{m}
\}\{ sa_1,\dots,sa_n\} =0\]
(see Sections~\ref{rem} and~\ref{kare} for the notation $\tilde{m}$ and
the suspension operator $s$, which decreases the super degree by 1);
one makes an $\Li$ algebra out of it via
\[ \lbr a_1,\dots,a_n\rbr =\{ m_n\} \{ a_1\} \cdots \{ a_n\} ,\]
as suggested in \cite{LS} and \cite{LM}. We will establish a master identity
\[ \{\tilde{m}\circ\tilde{m}\}\{ sa_1\}\{ sa_2\}\cdots\{ sa_n\} =0 \]
generating the usual $\Li$ identities for the higher brackets. 
Moreover, we will identify the Batalin-Vilkovisky
bracket for an odd linear operator $\tri$ and an even bilinear product $m_2$
as
\[ \{ a,b\}_{\tri}=(-1)^{|a|-1}[m_2,\tri ]\{ a,b\},\]
where $[m_2,\tri ]$ denotes the {\it Gerstenhaber bracket} defined on $\CA$ by
\[ {[x,y]}=x\circ y-(-1)^{d(x)d(y)+|x||y|}y\circ x.\]
As a result, the identities satisfied by $\{\; ,\;\}_{\tri}$ (and their
proofs) will be substantially simplified compared to~\cite{A}.
\vsp

Our goal is to apply these ideas eventually
to the ``homotopy'' structures on a TVOA
as in the work of Kimura, Voronov, and Zuckerman \cite{KVZ}.
For example, it is possible to go one step further and define ``partitioned
multilinear maps'' and their compositions, which will result in a brand new
master identity for homotopy Gerstenhaber algebras \cite{A2}. Another
-extremely difficult- project would be an algebraic construction of the
predicted higher products on the TVOA (precursors can be found in 
\cite{LZ}).  Note that although we stick to the complex 
number field throughout the article out of habit, all statements also hold for
fields of prime characteristic (except for cases with the
ubiquitous factor $1/2$), and algebraic closure is not required anywhere.

\section {The Hochschild complex of a graded vector space}

\subsection{Grading}\label{kare}

We will examine a context in which multilinear maps on a {\bf Z}-graded 
vector space
\[ A=\op_{j\in\Z}A^j\]
over {\bf C}, or more generally, linear maps
\[ x:TA\ra A\]
from the tensor algebra of $A$ into $A$ can be studied. We will assume that
the restriction $x_n$ of $x$ to $A^{\ot n}$ (not to be confused with the
homogeneous subspace $A^n$) is either homogeneous with respect to the (super)
{\bf Z}-grading, or else is a finite sum of homogeneous $n$-linear maps. The
notation for the {\bf super degree} will be
\[ |a|=j\;\;\;\mbox{if}\;\;\; a\in A^j,\]
and
\[ |x|=j\;\;\;\mbox{if}\;\;\; |x(a_1,\dots,a_n)|=|a_1|+\cdots +|a_n|+j\]
for all homogeneous $a_i\in A$ ($x:A^{\ot n}\ra A$). The terms ``odd
operator'' or ``even operator'' will refer to the super degree.
Most of the time we will reserve the name {\bf Hochschild complex} for
\be \CA =\Pi_{n=0}^{\infty}C^n(A)=Hom_{\C}(TA;A)=Hom_{\C}(\op_{n=0}^{\infty}
A^{\ot n};A) \label{one} \ee
instead of the classical
\be \CA =\op_{n=0}^{\infty}C^n(A)=\op_{n=0}^{\infty}Hom_{\C}(A^{\ot n};A),
\label{two}\ee
and will occasionally write formal expressions like
\[ x=x_1+x_2+\cdots \in\CA \]
(unfortunately, the subscripts will sometimes denote the corresponding
tensor power and sometimes an ordering of the symbols).
Note that (\ref{one}) can be thought of as a completion of (\ref{two}).  
In Section~\ref{twofour} we will allow
\[ \CBB =\HTT, \]
which has both (\ref{one}) and (\ref{two}) as subspaces, as well as $TA$.
\vsp

There is another natural concept of degree
on either type of ``cochains'' defined by
\[ D(x)=n\;\;\;\;\mbox{if $x$ is $n$-linear.}\]
Most of the time we will utilize the {\bf (adjusted) degree
of homogeneity}
\be d(x)=D(x)-1\label{three}\ee
instead for homogeneous elements of $\CA$, counting the number of tensor
factors in the domain of $x$ minus the number of tensor factors in the range.
If $R(x)$ denotes the tensor power of $A$ in the range of $x$,
the most general definition of $d(x)$ will be
\be d(x)=D(x)-R(x).\ee
At first glance, $d(x)$ appears only in
powers of $(-1)$ and can be replaced by $D(x)+1$ or even by one of $0$, $1$.
The particular choice (\ref{three}) will be most useful in what follows
below. We note that
\[ D(a)=0\;\;\;\;\mbox{and}\;\;\;\; d(a)=-1 \]
for $a\in A=C^0(A)$.
\vsp

Yet another degree associated with a graded vector space $A$ is the
so-called {\bf suspended (super) degree}
\be \| a\| =|a|-1.\label{yil}\ee
Under this shift, we denote the vector space by $sA$ and its elements $a$ by
$sa$ ($s$ is called the {\bf suspension operator}). In fact, $\| a\|$ is just
$|sa|$, if we think of $s:A\ra sA$ as a linear operator on $A$ with
\[ d(s)=0\;\;\;\mbox{and}\;\;\; |s|=-1\]
by some abuse of terminology. We will often omit the suspension operator
when ordinary round parentheses (as opposed to coupled, or curly,
parentheses) are used, as this notation does not involve hidden
permutations of symbols which would in turn necessitate sign changes. 
\vsp

Finally, some comments about terminology: the word ``super'' refers to the 
grading, while ``anti'' refers to the minus sign that is always present at
the exchange of two symbols. Then by ``super antisymmetry'' we mean
\[ ab=-(-1)^{|a||b|}ba.\]
We prefer to say a bilinear product on $A$
is super symmetric (or $|\;\;|$-graded symmetric) only when
\[ ab=(-1)^{|a||b|}ba.\]
Similarly, if $|\;\; |_1$ and $|\;\; |_2$ are two gradings on $A$, the identity
\[ ab=-(-1)^{|a|_1|b|_1+|a|_2|b|_2}ba\]
is bigraded antisymmetry, whereas
\[ ab=(-1)^{|a|_1|b|_1+|a|_2|b|_2}ba\]
is bigraded symmetry. To confuse matters further, note that the bigraded 
symmetry
\[ ab=(-1)^{|a||b|+d(a)d(b)}ba=-(-1)^{|a||b|}ba\]
of a bilinear product on $A$ is at the same time super antisymmetry.

\subsection{The coupled braces}

Most activities on algebras endowed with multilinear operations involve
compositions of maps and substitution of elements into maps. An excellent
notation of composition
\[  x \{ y\} =x\circ y\]
was invented by Gerstenhaber in \cite{Ger}, and
was generalized to 
\be  x \{ y_1,\dots,y_n\} \label{four}\ee
by Getzler in \cite{Get1}; we will prefer the uniform notation of 
{\bf coupled pairs of braces}
$\{ x\}\{ y\}$ and $\{ x\}\{ y_1,\dots,y_n\}$ advocated by Kimura, Voronov,
and Zuckerman in \cite{KVZ}. 
This last expression is a multilinear map
obtained by ``composing'' a multilinear map $x$ with similar maps
$y_1$,...,$y_n$ simultaneously (in this order). The idea is to generalize 
{\sl substitution}
\[ \{ x\} \{ a\} =x(a)\]
of $a\in A=C^0(A)$ into a linear map $x\in C^1(A)$, and {\sl composition}
\[ \{ x\}\{ y\} =x\circ y,\;\;\;\;\{ x\}\{ y\}\{ a\} =x(y(a))\]
of two linear operators $x$, $y\in C^1(A)$, paying attention to order and
grading. The resulting multilinear operator (\ref{four})  can be defined
again by ``graded and ordered substitution'' of elements of $A$, and in fact
the most complete expression involving multiple compositions/substitutions
will be of the form
\be \{ x\}\{ x_1^{(1)},\dots,x_{i_1}^{(1)}\}\cdots\{ x_1^{(m)},\dots,
x_{i_m}^{(m)}\}\{ a_1,\dots,a_n\}\in A, \label{five} \ee
where $x$, $x_j^{(i)}$ are (for the time being, bihomogeneous) elements of
$\CA$, possibly of $A$, and $a_1$,...,$a_n$ are enough (homogeneous)
elements of $A$ to fill the spaces allotted for arguments. In particular,
we are assuming that 
\[ R(x_1^{(1)})+\cdots +R(x_{i_1}^{(1)})=i_1\leq D(x)\] 
and so on, so that no symbols are left
out for lack of space at any stage of the substitution process (however,
see later expansion of this notation). The general
definition of (\ref{five}) is quite cumbersome but we can guess its form
from smaller examples. We first define
\be \{ x\}\{ a_1,\dots,a_n\}\stackrel{\rm def}{=}x(a_1,\dots,a_n)\in A
\label{six}\ee
for $x\in C^n(A)$, and put
\bea &\;\; &\{ x\}\{ x_1,\dots,x_m\}\{ a_1,\dots,a_n\} \label{seven}\\
&\stackrel{\rm def}{=}&\sum_{0\leq i_1\leq\cdots\leq i_m\leq n}
(-1)^{q(i_1,\dots,i_m)}
x(a_1,\dots,a_{i_1},x_1(a_{i_1+1},\dots),\dots,a_{i_m},x_m(a_{i_m+1},\dots),
\dots,a_n)\nn\\ 
&=&\sum_{0\leq i_1\leq\cdots\leq i_m\leq n}(-1)^{q(i_1,\dots,i_m)}
\{ x\}\{ a_1,\dots,a_{i_1},\{ x_1\}\{ a_{i_1+1},\dots\} ,\dots,a_{i_m},\nn\\
&&\;\;\;\;\;\;\;\;\;\;
\;\;\;\;\;\;\;\;\;\;\{ x_m\}\{ a_{i_m+1},\dots,\} ,\dots,a_n\} ,\nn\eea
where
\bea q(i_1,\dots,i_m)&=&\sum_{p=1}^md(x_p)(d(a_1)+\cdots +d(a_{i_p}))
+\sum_{p=1}^m|x_p|(|a_1|+\cdots +|a_{i_p}|)\nn\\
&=&-\sum_{p=1}^md(x_p)i_p+\sum_{p=1}^m|x_p|(|a_1|+\cdots +|a_{i_p}|)\nn\eea
denotes the sign change due to $x_p$ ``passing through'' $a_1$,...,$a_{i_p}$
in a deviation from the prescribed order on the left hand side. In both
(\ref{six}) and (\ref{seven}) the  $d$-grading of the ingredients add
up to that of the finished product (namely $-1$), and indeed even incomplete
expressions of coupled braces 
preserve the adjusted degree (\ref{three}) of homogeneity. In an expression
like (\ref{five}) consisting of homogeneous arguments, the number $n$ which
makes (\ref{five}) an element of $A$ can then be determined from
\[ d(x)+\sum_{i,j}d(x_j^{(i)})-n=-1.\]
We can now describe what (\ref{five}) ought to be by looking at (\ref{six})
and (\ref{seven}). Let us call the ordered elements of $\CA$ inside any
pair of braces a {\bf string}, and agree that 

(i) strings to the left contain ``higher'' entries than strings to the right,
and all entries in the same string are equivalent in ``height'' (exception:
all elements of $A$ have the same, and lowest, height);

(ii) every ``lower'' entry must appear in a ``higher'' entry (unless both 
entries are in $A$), not necessarily an adjacent one; and

(iii) order within any one string must be preserved.\newline
When we add up all possible expressions with the correct exchange signs,
we obtain the definition of (\ref{five}). It is possible to erase some
entries in a string, or even complete strings of elements of $A$ from the
right, while replacing the 
corresponding $d$- and super degrees by zero, and get a meaningful
multilinear operator. The ``higher pre-Jacobi identity''
\bea &\;\; &\{ x\}\{ x_1,\dots,x_m\}\{ y_1,\dots,y_n\} \label{eight}\\
&=&\sum_{0\leq i_1\leq\cdots\leq i_m\leq n}(-1)^{q(i_1,\dots,i_m)}
\{ x\}\{y_1,\dots,y_{i_1},\{ x_1\}\{ y_{i_1+1},\dots\} ,\dots,y_{i_m},\nn\\
&&\;\;\;\;\;\;\;\;\;\;\;\;\;\;\;\;\;\;\;\;
\{ x_m\}\{ y_{i_m+1},\dots\} ,\dots,y_n\}\nn\eea
of Voronov and Gerstenhaber \cite{VG} with
\[ q(i_1,\dots,i_m)=\sum_{p=1}^md(x_p)(d(y_1)+\cdots +d(y_{i_p}))+
\sum_{p=1}^m|x_p|(|y_1|+\cdots +|y_{i_p}|)\]
is the perfect example ($\{ z\}\{ \}$
means just $\{ z\}$). 
\vsp

We will also encounter expressions like
\[ \{ m\}\{ a,\;\;\} =\{ m\}\{ a,\mbox{id}\} ,\]
indicating, for example, that the first argument of a bilinear map $m$
is fixed. To this end, let us introduce the {\bf adjoint} of
an $n$-linear operator $x$, namely
\[ \adj(x):A^{\ot(n-1)}\ra C^1(A),\]
by defining
\be \{\adj(x)\{ a_1,\dots,a_{n-1}\}\}\{ a\} =x(a_1,\dots,a_{n-1},a).\ee
Note that the adjoint operator indicates bracketing by an element on the
left in case of a Lie algebra, with the Lie bracket as $x$.
\vsp

The definition of (\ref{five}) can be extended to nonhomogeneous
elements of the classical complex (\ref{two}) and even to those in 
(\ref{one}), as will be shown in Section~\ref{twofour}. We will
treat this definition as the mother lode of all definitions and identities
to come.

\begin{Rm} \label{firstrem} As is standard in algebraic topology,
the suspension notation affects the crossing of symbols in $C(sA)$;
the exchange of $sa$ and $sb$ is accompanied by 
\[ (-1)^{\| a\|\| b\|}=(-1)^{(|a|+1)(|b|+1)}\]
rather than 
\[ (-1)^{|a||b|+1} .\] 
For the rest of the exchange rules in $C(sA)$ see Section~\ref{rem}.
\end{Rm}

\subsection{Iterated braces}

The issue of {\bf iterated braces} must be regarded with caution. A grouping
of several strings within an extraneous pair of braces should simply mean
``make these substitutions first''. A well-known case
is the ``pre-Jacobi identity'' 
\be (x\circ y)\circ z-x\circ (y\circ z)=(-1)^{d(y)d(z)+|y||z|}((x\circ z)
\circ y-x\circ (z\circ y)),\label{nine}\ee
or
\be \{\{ x\}\{ y\}\}\{ z\} -\{ x\}\{\{ y\}\{ z\}\} =(-1)^{d(y)d(z)+|y||z|}
(\{\{ x\}\{ z\}\}\{ y\} -\{ x\}\{\{ z\}\{ y\}\})\label{vip}\ee 
in \cite{Ger}, for the {\bf Gerstenhaber product}
\be x\circ y=\{ x\}\{ y\}\label{ten}\ee
on $\CA$ (again defined in \cite{Ger}) between multilinear maps $x$ and
$y$. We note in passing that (\ref{nine}) is the defining identity
for a (bigraded) {\bf right pre-Lie algebra} $(B,\circ)$ as in \cite{Ger},
and simply says
\be {[R_y,R_z]}+R_{[y,z]}=0,\label{otf}\ee
where $R_y$ denotes right multiplication by $y$ in $B$, the brackets on the
left denote the usual super and $d$-graded commutator in 
$End(B)$ (we are thinking of $B=\CA$),
and the brackets on the right denote the {\bf Gerstenhaber bracket}
\be {[y,z]}\stackrel{\rm def}{=}y\circ z-(-1)^{d(y)d(z)+|y||z|}z\circ y.
\label{eleven}\ee
It is interesting that $\CA$ is {\sl not} a left pre-Lie algebra, i.e. the
identity
\[ {[L_y,L_z]}=L_{[y,z]}\]
involving left multiplications does not hold! 
That $\CA$ is a right pre-Lie algebra is proven in detail in \cite{Ger}.
Here is a more intuitive proof.

\begin{lemma}[Gerstenhaber \cite{Ger}]\label{prelie}
The Hochschild complex (\ref{one}) of a vector space is a right pre-Lie
algebra with respect to the G-product.
\end{lemma}

{\em Proof.} First of all we note that there is no point in enclosing
two strings from the {\sl left} in an extra pair of braces, as there is
only one possibility for substitution. The left
hand side of the identity~(\ref{nine}), namely
\[ \{ x\}\{ y\}\{ z\} -\{ x\}\{ \{ y\}\{ z\}\}\;\;\;
\mbox{(applied to some $\{ a_1,\dots,a_n\}$),}\]
consists of terms in which $z$ does {\sl not} appear inside $y$, i.e. in which
$y$ and $z$ appear in different entries of $x$. The right hand side, namely
\[ (-1)^{d(y)d(z)+|y||z|}(\{ x\}\{ z\}\{ y\} -\{ x\}\{\{ z\}\{ y\}\}),\]
consists of terms in which $y$ does not appear inside $z$, or again terms
for which $y$ and $z$ appear separately inside $x$. The sign rules are
the same on either side of~(\ref{nine}), and the sign on the right hand 
side takes care
of the initial misordering with respect to the left hand side.$\Box$ 
\vsp

We would like to generalize this statement to 
\[ \{ x\}\{ y\}\{ z_1,z_2\} -\{ x\}\{\{ y\}\{ z_1,z_2\}\} .\]
A quick inspection reveals that the left hand side contains 
extra terms in which exactly one of $z_1$ and $z_2$ appears inside $y$.
Keeping in mind the {\sl meaning} of coupled braces and the sign conventions,
we can easily concoct a valid identity:
\bea &&\{ x\}\{ y\}\{ z_1,z_2\} -\{ x\}\{\{ y\} ,\{ z_1,z_2\}\}\nn\\
&&-\{ x\}\{\{ y\}\{ z_1\} ,z_2\} -(-1)^{d(y)d(z_1)+|y||z_1|}\{ x\}\{ z_1,
\{ y\}\{ z_2\}\}\nn\\
&=&(-1)^{d(y)(d(z_1)+d(z_2))+|y|(|z_1|+|z_2|)}(\{ x\}\{ z_1,z_2\}\{ y\}
-\{ x\}\{\{ z_1,z_2\}\{ y\}\}),\nn\eea
or expanding the very last term,
\bea \label{tri}&&\{ x\}\{ y\}\{ z_1,z_2\} -\{ x\}\{\{ y\}\{ z_1,z_2\}\}\\
&&-\{ x\}\{\{ y\}\{ z_1\} ,z_2\} -(-1)^{d(y)d(z_1)+|y||z_1|}\{ x\}\{ z_1,
\{ y\}\{ z_2\}\}\nn\\
&=&(-1)^{d(y)(d(z_1)+d(z_2))+|y|(|z_1|+|z_2|)}(\{ x\}\{ z_1,z_2\}\{ y\}\nn\\
&&-(-1)^{d(y)d(z_2)+|y||z_2|}\{ x\}\{\{ z_1\}\{ y\} ,z_2\} -\{ x\}\{ z_1,
\{ z_2\}\{ y\}\} ).\nn\eea
Note that identity~(\ref{tri}) can also be put into the following form.

\begin{lemma} \label{triple} In $\CA$ we have
\bea &&\{ x\}\{ y\}\{ z_1,z_2\} -\{ x\}\{\{ y\}\{ z_1,z_2\}\}
-\{ x\}\{ [y,z_1],z_2\} -(-1)^{d(y)d(z_1)+|y||z_1|}\{ x\}\{ z_1,[y,z_2]\}\nn
\\&=&(-1)^{d(y)(d(z_1)+d(z_2))+|y|(|z_1|+|z_2|)}
\{ x\}\{ z_1,z_2\}\{ y\} \nn\eea
as an analogue of the pre-Jacobi identity.
\end{lemma}

By the way, identity~(\ref{nine}) has a similar presentation, namely

\begin{lemma} The pre-Jacobi identity can be written as
\[ \{ x\}\{ y\}\{ z\} -\{ x\} [y,z]
=(-1)^{d(y)d(z)+|y||z|}\{ x\}\{ z\}\{ y\} \]
for $x$, $y$, $z\in\CA$.
\end{lemma}

We can also give an elementary proof of

\begin{lemma}[Gerstenhaber \cite{Ger}]
A right (or left) bigraded pre-Lie algebra $(B,\circ)$ is a bigraded 
Lie algebra with respect to the Gerstenhaber bracket (\ref{eleven}).
\end{lemma}

{\em Proof.} The G-bracket is obviously bigraded antisymmetric. Moreover,
the bigraded cyclic Jacobi identity
\be (-1)^{d(x)d(z)+|x||z|}[[x,y],z]+(-1)^{d(y)d(x)+|y||x|}[[y,z],x]
+(-1)^{d(z)d(y)+|z||y|}[[z,x],y]=0\label{jac}\ee
is satisfied (again, the sign convention is that of $B=\CA$):
omitting the predictable signs due to the bigrading, we rewrite the left 
hand side of (\ref{jac}) as
\bea &&R_z([x,y])-R_{[x,y]}(z)+\mbox{cyclic}\nn\\
&=&R_z(R_y(x)-R_x(y))-R_{[x,y]}(z)+\mbox{cyclic}\nn\\
&=&(R_zR_y-R_yR_z-R_{[z,y]})(x)+\mbox{cyclic},\nn\eea
which is zero by the pre-Lie condition~(\ref{otf}). For another proof 
of the fact that $\CA$ is a graded Lie algebra see~\cite{Jim}.$\Box$
\vsp

One of the first places where 
the concept of a graded Lie algebra was introduced is Gerstenhaber's
\cite{Ger}. K.~Haring \cite{Ha} (who made some living history
investigations) clarifies the history of graded Lie algebras in her UNC
Master's Thesis. Another
source for the term and the abbreviation GLA is \cite{FN} by Fr\"{o}licher
and Nijenhuis.

\begin{ex}
As a demonstration of iterated braces, let us write down the first
two terms of (\ref{nine}) for say bilinear $x$, $y$, $z$ and elements $a$,
$b$, $c$, $d$ of $A$: the first one is
\bea &&\{\{ x\}\{ y\}\}\{ z\}\{ a,b,c,d\}
=\{ x\}\{ y\}\{ z\}\{ a,b,c,d\}\nn\\
&=&\{ x\}\{ y\}\{ z(a,b),c,d\} +(-1)^{-d(z)+|z||a|}\{ x\}\{ y\}\{ a,z(b,c),d\}
\nn\\ &&+(-1)^{-2d(z)+|z|(|a|+|b|)}\{ x\}\{ y\}\{ a,b,z(c,d)\}\nn\\
&=&\{ x\}\{ y(z(a,b),c),d\} +(-1)^{-d(y)+|y|(|z|+|a|+|b|)}\{ x\}\{ z(a,b),
y(c,d)\}\nn\\
&&+(-1)^{-d(z)+|z||a|}\{ x\}\{ y(a,z(b,c)),d\} +(-1)^{-d(z)+|z||a|-d(y)+|y||a|}
\{ x\}\{ a,y(z(b,c),d)\}\nn\\
&&+(-1)^{|z|(|a|+|b|)}\{ x\}\{ y(a,b),z(c,d)\} +(-1)^{|z|(|a|+|b|)-d(y)+|y||a|}
\{ x\}\{ a,y(b,z(c,d))\}\nn\\
&=&x(y(z(a,b),c),d)\pm x(z(a,b),y(c,d))\label{twelve}\\
&&\pm x(y(a,z(b,c)),d)\pm x(a,y(z(b,c),d))\nn\\
&&\pm x(y(a,b),z(c,d))\pm x(a,y(b,z(c,d))).\nn \eea
Meanwhile, a similar calculation
gives
\bea &&\{ x\}\{\{ y\}\{ z\}\}\{ a,b,c,d\}\nn\\
&=&\{ x\}\{\{ y\}\{ z\}\{ a,b,c\} ,d\}\pm\{ x\}\{ a,\{ y\}\{ z\}\{ b,c,d\}\}
\label{thirteen}\\
&=&x(y(z(a,b),c),d)\pm x(y(a,z(b,c)),d)\nn\\
&&\pm x(a,y(z(b,c),d))\pm x(a,y(b,z(c,d))).\nn\eea
\end{ex}

Iterated braces at two different levels -such as those on $A$ and $TA$- will 
be studied later.

\subsection{Antisymmetrization}

An unexpected bonus of the coupled-braces notation is the elimination of
explicit super antisymmetrizations. If $D(m)=n$, by definition we have
\be \{ m\}\{ a_1\}\cdots\{ a_n\} =\sum_{\sigma\in S_n}(-1)^{p(\sigma ;a_1,
\dots,a_n)}
m(a_{\sigma(1)},\dots,a_{\sigma(n)}),\label{fifteen}\ee
where $S_n$ is the symmetric group, and
\bea p(\sigma ;a_1,\dots,a_n)&\defi&
\sum_{u<v,\sigma^{-1}(u)>\sigma^{-1}(v)}d(a_u)d(a_v)+
\sum_{u<v,\sigma^{-1}(u)>\sigma^{-1}(v)}|a_u||a_v|\nn\\
&=&\mbox{$\#$ of exchanges}+e(\sigma ;a_1,\dots,a_n),\nn\eea
so that
\[ (-1)^{p(\sigma ;a_1,\dots,a_n)}\;\defi\; {\rm sgn}(\sigma){\ep(\sigma ;
a_1,\dots,a_n)},\]
where
\[ e(\sigma ;a_1,\dots,a_n)\,\defi\,\sum_{u<v,\sigma^{-1}(u)>\sigma^{-1}(v)}
|a_u||a_v| \]
and
\[\epsilon(\sigma ;a_1,\dots,a_n)\,\defi\, (-1)^{e(\sigma ;a_1,\dots,a_n)}.\]
For the suspended grading, we have the analogous definition
\bea \tep(\sigma ;a_1,\dots,a_n)&\defi&\sum_{\uv}\| a_u\|\;\| a_v\|\\
&=&\sum_{\uv}(|a_u|-1)(|a_v|-1)\nn\\
&=&p(\sigma ;a_1,\dots,a_n)+\sum_{\uv}(|a_u|+|a_v|)\;\;\;\mbox{(mod 2)}\nn\\
&=&p(\sigma ;a_1,\dots,a_n)+\sum_{\uv}(\| a_u\| +\| a_v\|)\;\;\;
\mbox{(mod 2)},\nn\eea
with
\be \tilde{\epsilon}(\sigma ;a_1,
\dots,a_n) \;\defi\; (-1)^{\tep(\sigma ;a_1,\dots,a_n)}.\ee
Then $p$ is good for antisymmetrizing a product
with respect to the super grading, and $\tep$ is good for symmetrizing 
with respect to the suspended grading! 
Note that
\be\label{split} p(\sigma\sigma';a_1,\dots,a_n)=
p(\sigma';a_1,\dots,a_n)+p(\sigma;a_{\sigma'(1)},\dots,
a_{\sigma'(n)})\;\;\;\mbox{(mod 2),}\ee
as we want to add up exchange terms coming from $\sigma'$ followed by more
exchanges coming from $\sigma$, and whenever $\sigma$ unravels an exchange
done by $\sigma'$, the subtotal is zero modulo~2. The same goes for $e$
and $\tep$, and it is well-known that ${\rm sgn}(\sigma\sigma')={\rm sgn}(
\sigma){\rm sgn}(\sigma')$. 
As a special case of antisymmetrization,
\be \lbr a,b\rbr\stackrel{\rm def}{=}
\{ m\}\{ a\}\{ b\} =m(a,b)-(-1)^{|a||b|}m(b,a)\label{sixteen}\ee
is the usual definition of a bilinear bracket associated to a bilinear
product $m$, satisfying the super Jacobi identity when $m$ is associative. In 
general, (\ref{fifteen}) defines a $|\;\; |$-graded antisymmetric product
\be \lbr a_1,\dots,a_n\rbr\stackrel{\rm def}{=}\{ m\}\{ a_1\}\cdots\{ a_n\}
\label{seventeen}\ee
on $A$ with
\be \lbr a_{\sigma(1)},\dots,a_{\sigma(n)}\rbr =(-1)^{p(\sigma ;a_1,\dots,
a_n)}\lbr a_1,\dots,a_n\rbr ={\rm sgn}(\sigma)\epsilon(\sigma ; 
a_1,\dots,a_n)\lbr a_1,\dots,a_n\rbr\label{ds}\ee
for any $\sigma\in S_n$.
We will come back to this bracket in connection with $\Ai$ and $\Li$
algebras. For later use, we have

\begin{lemma}\label{lu} The expression
\be (-1)^{\su +\pos} \label{expr}\ee
is independent of $\sigma$. In particular, it is equal to
\[ (-1)^{\po}.\]
\end{lemma}

{\it Proof.} Appendix.$\Box$

\subsection{Derivations and higher order differential operators}

\label{phi}
An inductive definition of higher order differential operators on a 
superalgebra $A$ with a (noncommutative, nonassociative) bilinear product
$m=m_2$, consistent with the commutative and associative case described by
Koszul in \cite{Ko}, was given in \cite{A} and was shown to be a good fit
for modes of vertex operators. Namely, a (homogeneous)
linear operator 
$\tri :A\ra A$ is a {\bf differential operator of order r} and of super degree
$|\tri |$ if and only if
\[ \Phi_{\tri}^{r+1}(a_1,\dots,a_{r+1})=0\;\;\;\;\forall a_i\in A,\]
where
\bea &&\Phi_{\tri}^1(a)=\tri(a),\nn\\
&&\Phi_{\tri}^2(a,b)=\Phi_{\tri}^1(ab)-\Phi_{\tri}^1(a)b-(-1)^{|a||\tri |}
a\Phi_{\tri}^1(b),\nn\\
&&\vdots\nn\\
&&\Phi_{\tri}^{r+1}(a_1,\dots,a_{r+1})=\Phi_{\tri}^{r}(a_1,\dots,a_ra_{r+1})
-\Phi_{\tri}^r(a_1,\dots,a_r)a_{r+1},\nn\\
&&-(-1)^{|a_r|(|\tri |+|a_1|+\cdots +|a_{r-1}|)}a_r\Phi_{\tri}^r(a_1,\dots,
a_{r-1},a_{r+1})\nn\\
&&\vdots\nn\eea
($m$ suppressed). The multilinear forms $\Phi_{\tri}^r$ can be expressed
as follows in the coupled-braces notation:
\bea \Phi_{\tri}^1(a) &=&\{\tri\}\{ a\} ,\nn\\
\Phi_{\tri}^2(a,b) &=&[\Phi_{\tri}^1,m_2]\{ a,b\} ,\;\;\;\mbox{and}\nn\\
\Phi_{\tri}^{r+2}(a_1,\dots,a_r,a,b)
&=&{[\{\Phi_{\tri}^{r+1}\}\{ a_1,\dots,a_r,\mbox{id}\} ,m_2]}\{ a,b\}\nn\\
&=&\{\Phi_{\{\Phi_{\tri}^{r+1}\}\{ a_1,\dots,a_r,{\rm id}\}}^2\}\{ a,b\}
\;\;\;\mbox{for $r\geq 1$.}\nn\eea 
Alternatively, in terms of the adjoint operators, we have
\bea \Phi_{\tri}^{r+2}(a_1,\dots,a_r,a,b) &=& [\adj(\Phi_{\tri}^{r+1})
\{ a_1,\dots,a_r\} ,m_2]\{ a,b\}\label{ooo}\\
&=& \{\Phi^2_{\adj(\Phi_{\tri}^{r+1})\{ a_1,\dots,a_r\} }\}\{ a,b\} .\nn\eea
In particular, the linear operator $\tri$ is a derivation
of $m_2$ if and only if the Gerstenhaber bracket ${[\tri,m_2]}$ is
identically zero. 
\vsp

\begin{lemma}\label{yasa} For odd linear operators $T$ and $U$ on $A$, the 
bracket $[T,U]=TU+UT$ 
is related to the Gerstenhaber brackets of the $\Phi$ operators
as follows:
\bea \Phi_{[T,U]}^1(a)&=&[\Phi_T^1,\Phi_U^1]\{ a\}\nn\\
\Phi_{[T,U]}^2(a,b)&=&[\Phi_T^1,\Phi_U^2]\{ a,b\} +[\Phi_U^1,\Phi_T^2]
\{ a,b\}\nn\\
\Phi_{[T,U]}^3(a,b,c)&=&[\Phi_T^1,\Phi_U^3]\{ a,b,c\} +[\Phi_U^1,\Phi_T^3]
\{ a,b,c\}
\nn\\ &&+[\Phi_T^2,\adj(\Phi_U^2)\{ a\} ]\{ b,c\} +[\Phi_U^2,\adj(\Phi_T^2)\{ a
\} ]\{ b,c\} .\nn\eea
\end{lemma}

{\it Proof.} Straightforward.$\Box$
\vsp

With the new coupled braces, it is easy to generalize the idea of higher
order differential operators $\tri$ with respect to a bilinear map $m_2$
to {\sl multilinear maps} which are differential operators with respect to
another multilinear map! The obvious way is to introduce new operators
\be \oper{m_l}{r}{m_k} \label{ope}\ee
where $m_k$ and $m_l$ are $k$-linear and $l$-linear maps respectively, and
$r$ is once again a positive integer. When $l=2$ and $m_k=\tri$ is a linear
map, (\ref{ope}) will coincide with $\Phi_{\tri}^r$.
We make the inductive definition
\bea \oper{m_l}{1}{m_k}(a_1,\dots,a_k) &=&\{ m_k\}\{ a_1,\dots,a_k\} ,\\
\oper{m_l}{2}{m_k}(a_1,\dots,a_{k+l-1}) &=&[m_k,m_l]\{ a_1,\dots,
a_{k+l-1}\} ,\nn\eea
and
\bea &&\oper{m_l}{r+2}{m_k}(a_1,\dots,a_{(r+1)(l-1)+k})\label{oa} \\
&=&[\adj(\oper{m_l}{r+1}{m_k})
\{ a_1,\dots,a_{r(l-1)+k-1}\} ,m_l]\{ a_{r(l-1)+k},\dots,
a_{(r+1)(l-1)+k}\} \nn\\
&=&\{\oper{m_l}{2}{\adj(\oper{m_l}{r+1}{m_k})\{ a_1,\dots,a_{r(l-1)+k-1}\} }\}
\{ a_{r(l-1)+k},\dots,a_{(r+1)(l-1)+k}\}
\;\;\;\mbox{for $r\geq 1$}.\nn\eea
Note that
\be d(\oper{m_l}{r}{m_k})=(r-1)d(m_l)+d(m_k)=(r-1)(l-1)+k-1,\ee
and
\be |\oper{m_l}{r}{m_k}|=(r-1)|m_l|+|m_k|.\ee
Predictably, we define higher order multilinear differential operators by

\begin{defn} A $k$-linear map $m_k$ is a {\bf differential operator of
order $r$} with respect to an $l$-linear map $m_l$ if and only if
$\oper{m_l}{r+1}{m_k}$ is identically zero.
\end{defn}

\begin{Rm} The operator (\ref{ope}) is linear in $m_k$. It is not symmetric
in $m_k$ and $m_l$ (except for $r=2$) and it is definitely
biased, because of the lopsided adjoint operator.
Moreover, the definition would improve if we
write $\oper{m_k}{r}{m_l}$ and reverse the arguments of the G-brackets 
in~(\ref{oa}), for then we would have the exact same ordering of symbols
on both sides of the definition. Nevertheless, this version would differ
from~(\ref{ooo}) and (\ref{oa}) only by an overall minus sign (provided
$m_l$ is even in the second case).
\end{Rm}

\subsection{Modified multilinear maps}

\label{rem} We will encounter many examples of modification $\tilde{m}$ of
a multilinear map $m$ by a sign that depends on the grading of the arguments. 
Roughly speaking, this modification translates between multilinear maps on 
graded symmetric and graded exterior algebras on the same underlying vector
space $A$ with two different gradings ($\|\;\|$ goes with symmetric and $|\;|$
goes with antisymmetric). More precisely, we expect one multilinear map
(say $m$),
even if not antisymmetric itself, to satisfy some identities in which an
exchange of $a$ and $b$ is accompanied by
\[ (-1)^{|a||b|+d(a)d(b)}=-(-1)^{|a||b|}\]
(we may also say these identities are ``bigraded'', in the sense of super and
$d$-gradings). Meanwhile, $\tilde{m}$ will satisfy a similar identity
in which the exchange of $sa$ and $sb$ will be marked by the factor
\[ (-1)^{\| a\|\;\| b\|}.\]
(Kjeseth's thesis \cite{Kje} 
and Penkava's article \cite{Pen1} carefully explain the interplay between
the symmetric and antisymmetric settings, or between $C(A)$ and $C(sA)$.)
The exact factor of modification from $m$ to $\tilde{m}$ was most clearly 
stated 
in \cite{GS} (Lemma~1.3). We can see why it is required from the following 
coupled-braces argument:
if $s:A\ra sA$ is the suspension operator $a\mapsto sa$, we
have to define the isomorphism
\[ s^{\ot n}:A^{\ot n}\ra(sA)^{\ot n}\]
by
\bea &&s^{\ot n}(a_1\ot a_2\ot\cdots\ot a_n)\nn\\
&=& \{ s,\dots,s\}\{ a_1,\dots,a_n\}\;\;\;\mbox{(see
Section~\ref{twofour} for expanded notation)}\nn\\
&=&(-1)^{|a_1| +(| a_1| +| a_2|)+\cdots+(| a_1| +| a_2| +\cdots
+| a_{n-1}|)}sa_1\ot\sa_2\ot\cdots\ot sa_n\nn\\
&=&(-1)^{\sum_{t=1}^n(n-t)| a_t|}sa_1\ot\cdots\ot sa_n.\nn\eea
Now if
\[ m:A^{\ot n}\ra A\]
is an $n$-linear operator, we define its counterpart
\[ \tilde{m}:(sA)^{\ot n}\ra sA\]
by
\be \tilde{m}=s\circ\m\circ(s^{-1})^{\ot n},\label{tildem}\ee
or by
\bea\label{til} \tilde{m}(a_1,\dots,a_n)
&=&\{\tilde{m}\}\{ sa_1,\dots,sa_n\}\\
&=&\{ s\}\{ m\}\{ s^{-1},\dots,s^{-1}\}\{ sa_1,\dots,sa_n\}\nn\\
&=&(-1)^{| s^{-1}|\sum_{t=1}^n|sa_t|}\{ s\}\{ m\}\{ a_1,\dots,a_n\}\nn\\
&=& (-1)^{\sum_{t=1}^n(n-t)\|a_t\|}m(a_1,\dots,a_n).\nn\eea

\vsp

We would like to determine the exchange rules among symbols like $sa$
and $\tm$ in $C^{\bullet}(sA)$ in accordance with the old rules. We claim
that replacing the bidegree with the {\bf suspended degrees}
\be \| a\| =|a|-1\;\;\;\mbox{and}\;\;\; \| m\| =|m|+d(m)\ee
of $sa$ and $\tm$
is sufficient. Note that since both $d$ and the super degree are preserved
by the coupled braces, so is the grading $\|\;\;\|$. We will not give a
complete proof of the correctness of this translation, but rather provide 
individual cases of justification. Of course, one may also adopt these
exchange rules as the {\sl definition}. For example, we expect

\begin{lemma} $\{\tm\}\{ sa\}\{ sb\} =(-1)^{\| a\|\;\| b\|}\{\tm\}\{ sb\}
\{ sa\}$\end{lemma}

for bilinear $m$. Indeed,
\vsp

{\it Proof.} 
\bea \{\tm\}\{ sa\}\{ sb\} &=& (-1)^{|a|}\{\tm\}\{ s,s\}\{ a\}\{ b\}\nn\\
&=&(-1)^{|a|+|a||b|+1}\{\tm\}\{ s,s\}\{ b\}\{ a\}\nn\\
&=&(-1)^{|a|+|b|+|a||b|+1}(-1)^{|b|}\{\tm\}\{ s,s\}\{ b\}\{ a\}\nn\\
&=&(-1)^{\| a\|\;\| b\|}\{\tm\}\{ sb\}\{ sa\} .\Box\nn\eea

Hence
\[ \{\tilde{m}\}\{ sa_1\}\cdots\{ sa_n\}\]
is a {\it symmetric} product. Note that

\begin{lemma} $\{\tm\}\{ sa_1\}\cdots\{ sa_n\} =(-1)^{\sum_{t=1}^n(n-t)
\| a_t\|}\{ s\}\{ m\}\{ a_1\}\cdots\{ a_n\} .$\end{lemma}

{\it Proof.} Omitted. It is similar to the proof of Proposition~\ref{pro8}
on strongly homotopy Lie algebras.$\Box$

\begin{lemma}\label{supsym} If $m$ is super antisymmetric, then $\tilde{m}$ is 
suspended-graded symmetric.\label{mti}
\end{lemma}

{\it Proof.} Appendix.$\Box$
\vsp

It is also best to define the Gerstenhaber bracket of $\tm$ and $\tx$ as
\be {[\tm,\tx ]}\;\defi\;\tm\circ\tx -(-1)^{\| m\|\;\| x\|}\tx\circ\tm,\ee
because then we have -among other consistency conditions-

\begin{lemma} $[\tm,\tx ]\{ sa_1,\dots,sa_n\} =(-1)^{d(m)\| x\| +\sum_{t=1}^n
(n-t)\| a_t\|}\{ s\} [m,x]\{ a_1,\dots,a_n\}$.\label{comut}\end{lemma}

{\it Proof.} Appendix.$\Box$

\subsection{Extension of coupled braces} 
\label{twofour}

\subsubsection{Derivations and coderivations of the tensor algebra}

We would like to venture beyond the conventional use of coupled braces for
multilinear maps with values in $A$, and expand our notation to handle
derivations and coderivations of 
\[ T^{\bullet}A=\op_{n=0}^{\infty}A^{\ot n}.\]
The multiplication on $TA$ is given by
\[ M\in Hom(TA\ot TA;TA),\]
with
\bea M(a_1\ot\cdots\ot a_k,a_{k+1}\ot\cdots\ot a_n)&=&\{ M\}'\{\{
a_1,\dots,a_k\},\{ a_{k+1},\dots,a_n\}\}'\nn\\
&\defi&\{\{ a_1,\dots,a_k,a_{k+1},\dots,a_n\}\}' \nn\eea
(primed braces live in $T(TA)$ as opposed to $TA$); the symbol
\[ \{ a_1,\dots,a_n\}\]
by itself was obviously meant to be
\[ a_1\ot\cdots\ot a_n\in A^{\ot n}\]
all along. Meanwhile
\[ \{ a_1,\dots,a_k\}\{ a_{k+1},\dots,a_n\}\]
is a signed sum in $A^{\ot n}$ over all tensor products of the $a_i$'s
preserving the order in both strings -also called {\it shuffles}
(recall that all $a_i$ are of the same
``height'' and cannot be substituted into each other),
and the {\bf second level braces} $\{\;,\;
\}'$ belong to tensor products in $T(TA)$ and operators with values in $T(TA)$.
We will allow multilinear maps to take values in $TA$, and define
\[ d(a_1\ot\cdots\ot a_n)=-n\;\;\;\mbox{for}\;\;\; a_i\in A,\]
and
\[ d(x)=D(x)-R(x)=k-l\;\;\;\mbox{for}\;\;\; x:A^{\ot k}\ra A^{\ot l},\]
perfectly consistent with our earlier conventions (coupled braces still 
preserve $d$). 

\begin{Rm} The pre-Jacobi identity (\ref{vip}) holds for $x=a$, $y=b$, and
$z=c$, as the product $a\circ b=\{ a\}\{ b\}$ is associative on $A$; both
sides of the identity vanish. It is also easily checked that the
Gerstenhaber bracket is identically zero on $A\ot A$. This is again entirely
consistent with the old complex, where the grading $d=-2$ does not exist.
\end{Rm}

The comultiplication on $TA$ is given by the diagonal map
\[ \De\in Hom(TA;TA\ot TA)\]
with
\bea \De(a_1\ot\cdots\ot a_n)&=&\{\De\}'\{\{ a_1,\dots,a_n\}\}'\nn\\
&\defi&\sum_{k=0}^n(a_1\ot\cdots\ot a_k)\ot(a_{k+1}\ot\cdots\ot a_n)\nn\\
&=&\sum_{k=0}^n\{\{ a_1,\dots,a_k\},\{ a_{k+1},\dots,a_n\}\}' .\nn\eea

\begin{Rm} We will take $\{\;\} =1\in${\bf C}. Although some authors choose
to ignore the 0-th tensor power of $A$ in the context of higher homotopies,
we would like to include it for
completeness, as $A$ resides in $\CBB$ in the form $Hom(\mbox{{\bf C}};A)$.
\end{Rm}

Similar formulas hold for $sa_i\in sA$ if we replace $TA$ by $T(sA)$. It is
well-known that a derivation of $TA$ is determined by its restriction to
$A$, and
\[ Der(TA)\cong Hom(A;TA).\]
On the other hand, a coderivation of $TA$ is determined by itself followed
by the projection of $TA$ onto $A$, and
\[ Coder(TA)\cong Hom(TA;A)\]
(see Stasheff \cite{Jim}). Note that derivations ${\cal D}$ of $TA$ satisfy
\be {[{\cal D},M]}'=0,\label{derivation}\ee
and coderivations ${\cal C}$ of $TA$ satisfy
\be {[{\cal C},\De ]}'=0;\label{coderivation}\ee
the meaning of -second level- composition (and hence of Gerstenhaber 
bracket) between operators with values in the tensor algebra $T(TA)$ 
will become clear below. The primes refer to the fact that our
underlying vector space is not $A$ but $B=TA$.
\vsp

Let us define the {\bf extended Hochschild complex} by 
\be \CBB =\HTT ,\ee
and allow the arguments of coupled braces to live in $\CBB$. 
If $x:A\ra A^{\ot k}$ is a linear map,
we want to denote its extension to $TA$ as a derivation by
\bea x(a_1,\dots,a_n)&=&\{ x\}\{ a_1,\dots,a_n\}\\
&\defi&\sum_{i=1}^n(-1)^{d(x)(d(a_1)+\cdots +d(a_{i-1}))+|x|(|a_1|+\cdots +
|a_{i-1}|)}\{ a_1,\dots,a_{i-1},\{ x\}\{ a_i\},a_{i+1},\dots,a_n\}\nn\\
&=&\sum_{i=1}^n(-1)^{-d(x)(i-1)+|x|(|a_1|+\cdots +
|a_{i-1}|)}a_1\ot\cdots\ot x(a_i)\ot\cdots\ot a_n\nn\eea
on $A^{\ot n}$. This way, we allow the $a_i$'s to spread out to the left
and fill out {\sl every} available space, instead of defining the expression
$\{ x\}\{ a_1,\dots,a_n\}$ as zero.
The above formula has the generalization
\bea x(a_1,\dots,a_n)&=&
\{ x\}\{ a_1,\dots,a_n\}\\&\stackrel{\rm def}{=}&\sum_{i=1}^{n-r+1}
(-1)^{-d(x)(i-1)+|x|(|a_1|+\cdots +|a_{i-1}|)}a_1\ot\cdots\ot x(a_i,\dots,
a_{i+r-1})\ot\cdots\ot a_n\nn\eea
for $x:A^{\ot r}\ra A^{\ot k}$, $r<n$ (but the extension is not a derivation).
In this spirit, we find it natural to define
\bea \{ x,y\}\{ a_1,\dots,a_n\} &\defi&\pm\{\{ x\}\{ a_1,\dots,a_k\},
\{ y\}\{ a_{k+1},\dots,a_n\}\}\nn\\
&=&\pm x(a_1,\dots,a_k)\ot y(a_{k+1},\dots,a_n)\in TA\nn\eea
for $a_i\in A$, $x:A^{\ot k}\ra TA$, and $y:A^{\ot (n-k)}\ra TA$, and
set
\[ \{ x,y\}\{ a_1,\cdots,a_n\}\stackrel{\rm def}{=}
\sum_{i=1}^{n-1}\{ x_i,y_{n-i}\}\{ a_1,\dots,a_n\}\in TA\]
for $a_i\in A$, $x$, $y:TA\ra TA$. Note that we have already been using the
notation $\{ s,\dots,s\}$ to denote a signed tensor product of maps.
\vsp

In general, we define
\[ \{ x_1,\dots,x_n\}\{ y\} =\sum_{i=1}^n(-1)^{d(y)(d(x_{i+1})+\cdots +
d(x_n))+|y|(|x_{i+1}|+\cdots +|x_n|)}\{ x_1,\dots,\{ x_i\}\{ y\},\dots,
x_n\},\]
and
\bea &&\{ x_1,\dots,x_n\}\{ y_1,\dots,y_m\}\\ 
&=&\sum_{1\leq t\leq n;\; 1\leq i_1<\cdots <i_t\leq n;\; 1\leq j_1<\cdots 
<j_t=m}\nn\\
&&\;\;\;\;\;\pm\{ x_1,\dots,\{ x_{i_1}\}\{ y_1,\dots,y_{j_1}\},\dots,
\{ x_{i_2}\}\{ y_{j_1+1},\dots,y_{j_2}\},\dots,\{ x_{i_t}\}
\{ y_{j_{t-1}+1},\dots,y_{j_t}\},\dots,x_n\} .\nn\eea
\vsp

To summarize, the notation $\{ x\}\{ a_1,\dots,a_n\}$ for $x:A^{\ot k}\ra
A^{\ot l}$ consistently covers a variety of cases: if $k=n$ and $l=1$ (or even
$l>1$), these are exactly our old braces. If $n>k$, we let $x$ ``slide''
through the tensor product $a_1\ot\cdots\ot a_n$. If $x$ is a (finite or
infinite) sum of homogeneous parts $x_1$, $x_2$, ..., we write $x=x_1+x_2+
\cdots +x_n+\cdots$ 
and evaluate each piece accordingly. And as before, if $D(x)=k$
with $k>n$, we may continue to interpret $\{ x\}
\{ a_1,\dots,a_n\}$ as a $(k-n)$-linear map waiting to be fed (we sum over
all possible positions of $a_i$'s inside $x$ preserving the order),
or restrict the range and define it to be zero, depending on context.
With this summary we realize that we have, in fact, already crossed over
into the realm of coderivations! We define an extension of
\[ x:TA\ra A,\;\;\; x=x_1+x_2+\cdots\]
to
\[ x:TA\ra TA\]
as a coderivation by the well-known construction (again see \cite{Jim}), 
namely by
\bea x(a_1,\dots,a_n)&=&\{ x\}\{ a_1,\dots,a_n\}\nn\\
&=&\sum_{k=1}^n\{ x_k\}\{ a_1,\dots,a_n\} .\nn\eea

\subsubsection{Extension of the Gerstenhaber bracket}

Since the coupled braces $\{ x\}\{ y\}$ make sense for any two multilinear
maps $x$ and $y$ in $\CBB$, we may again define the G-bracket to be
\be {[x,y]}\;\defi\; x\circ y-(-1)^{d(x)d(y)+|x||y|}y\circ x.\ee
Let us look into the ``composition'' $\{ x\}\{ y\}$ more closely. If $R(y)\leq
D(x)$, then this expression is easy to figure out, as in
\[ \{ x\}\{ y\}\{ a,b,c,d\} =x(y(a,b),c,d)\pm x(a,y(b,c),d)\pm
x(a,b,y(c,d))\]
for
\[ x:A^{\ot 5}\ra A,\;\;\; y:A^{\ot 2}\ra A^{\ot 3},\;\;\; d(x)=4,\;\;\;
d(y)=-1.\]
If $R(y)>D(x)$, or the range of $y$ does not fit into the domain of $x$,
then
\[ \{ x\}\{ y\}\{ a_1,\dots,a_n\} \]
(say with $D(y)=n$) will be
\[ \{ x\}\{ y(a_1,\dots,a_n)\} \]
where again $x$ will slide over the tensors in $y(a_1,\dots,
a_n)$. In fact, we have the precise expressions
\be D(\{ x\}\{ y\})=D(y)+{\rm max}\{ 0,D(x)-R(y)\} \ee
and
\be R(\{ x\}\{ y\})=R(x)+{\rm max}\{ 0,R(y)-D(x)\} ,\ee
where -as before- $D$ denotes the exact number of arguments that would fill
every slot in the multilinear, tensor-valued map, and $R$ denotes the
tensor power in the range. This is still consistent with the ideas that
$d=D-R$ and $d$ is additive, because
\[ d(\{ x\}\{ y\})= D(\{ x\}\{ y\})-R(\{ x\}\{ y\})=D(y)+D(x)-R(y)-R(x)
=d(x)+d(y)\]
whether $D(x)>R(y)$ or not. Note that 
\[ D(\{ a\}\{ y\})=D(y)+{\rm max}\{ 0,-R(y)\} =D(y)\]
and
\[ R(\{ a\}\{ y\})=1+{\rm max}\{ 0,R(y)\} =R(y)+1,\]
and $\{ a\}\{ y\}$ means
\[ \{ a\}\{ y\}\{ a_1,\dots,a_n\} =\{ a\}\{ y(a_1,\dots,a_n)\} ,\]
a signed sum, for $D(y)=n$. For comparison, we note
\[ D(\{ y\}\{ a\})=D(y)-1\;\;\;\mbox{and}\;\;\; R(\{ y\}\{ a\})=R(y).\]
This discrepancy in domains did not exist in the old complex, as we always had
$R=1$, and
\[ D(\{ x\}\{ y\})=D(y)+D(x)-1=D(\{ y\}\{ x\}).\]
In the definition of the linear function $\de(a)$ in the following section,
we will only have
\[ {[m,a]}=\{ m\}\{ a\}\]
for some fixed $m$ with $D(m)=2$ and $R(m)=1$, because
\[ D(\{ a\}\{ m\})=R(\{ a\}\{ m\})=2,\]
hence $\{ a\}\{ m\}\not\in\CA$. Then
\[ \de(a)(b)=\{ m\}\{ a\}\{ b\} =ab-(-1)^{|a||b|}ba\]
is the correct formula (note that there will be two sign conventions for
the Hochschild cohomology differential $\de$).
Also recall that we have
\[ \{ a\}\{ b\} =\{ a,b\} -(-1)^{|a||b|}\{ b,a\} ,\]
with
\bea {[a,b]}&=&\{ a\}\{ b\} -(-1)^{|a||b|}\{ b\}\{ a\}\nn\\
&=&\{ a,b\} -(-1)^{|a||b|}\{ b,a\} +(-1)^{|a||b|}\{ b,a\} -\{ a,b\}\nn\\
&=& 0.\nn\eea
As for the pre-Jacobi identity for the new composition rule, 
the proof of Lemma~\ref{prelie} is valid word for word.
\vsp

With the extended definition of composition in mind, we can see why identities
(\ref{derivation}) and (\ref{coderivation}) are equivalent to the common
definitions of derivation and coderivation: the first one has already been 
commented on in Section~\ref{phi} on higher order differential operators. 
Identity~(\ref{coderivation}), on the other hand, is exactly 
\[ ({\cal C}\ot\ide +\ide\ot{\cal C})\circ\De =\pm\De\circ{\cal C},\]
the usual definition in the absence of our notation.

\section{Identities in various types of algebras}

\subsection{Associative algebras}
\label{onetwo}

For identities on an associative (super, or {\bf Z}-graded) algebra
$(A,m)$, with $m:A\ot A\ra A$, we can stick to the classical Hochschild
complex (\ref{two}) with coefficients in the two-sided module $A$.
We note that the associativity condition on $m$ can be written as
\be m\circ m=0.\label{twentytwo}\ee
Indeed,
\bea \{ m\}\{ m\}\{ a,b,c\}&\stackrel{\rm def}{=}&m(m(a,b),c)+(-1)^{d(m)
d(a)+|m||a|}m(a,m(b,c))\nn\\
&=&(ab)c-a(bc)\;\;\;\;\forall a,b,c\in A,\label{twentythree}\eea
as $d(m)=1$, $d(a)=-1$, and $|m|=0$ (this last one is an implicit assumption
in a superalgebra). It is also true that 
\[ \tilde{m}\circ\tilde{m}=0;\]
see (\ref{thirtyfive}).

\subsubsection{Classical definitions of the differential and the dot product}

Hochschild constructed a square-zero differential
\[ \de :C^n(A)\ra C^{n+1}(A)\]
on $\CA$ given by the formula
\bea (\de(x))(a_1,\dots,a_{n+1})&=&(-1)^{|a_1||x|}a_1x(a_2,\dots,a_{n+1})
\label{eighteen}\\
&&-x(a_1a_2,a_3,\dots,a_{n+1})+\cdots +(-1)^nx(a_1,a_2,\dots,a_na_{n+1})\nn\\
&&-(-1)^{n+1}x(a_1,\dots,a_n)a_{n+1},\nn\eea
where $m$ is suppressed. We implicitly understand that $a_1\in A$ is 
homogeneous and $x\in C^n(A)$ is bihomogeneous. Extension of the definition
to nonhomogeneous $x\in\CA$ and $a_1\in A$ is by linearity. Note that
\[ D(\de(x))=D(x)+1\;\;\;\;\mbox{and}\;\;\;\; |\de(x)|=|x|.\]
Identity (\ref{eighteen}) for $x=a\in C^0(A)=A$ is
\be \de(a)(b)=(-1)^{|a||b|}ba-ab=-(ab-(-1)^{|a||b|}ba).\label{nineteen}\ee
Then the algebra $(A,m)$ is super commutative (not super anticommutative!)
if and only if $\de :C^0(A)\ra C^1(A)$ is identically zero. Moreover, we have
\bea \de(x)(a,b)&=&(-1)^{|a||x|}ax(b)-x(ab)+x(a)b\label{twenty}\\
&=&-\Phi_x^2(a,b)\nn\eea
for linear $x$, and
\bea \de^2(a)(b,c)&=&(-1)^{|a|(|b|+|c|)}(b(ca)-(bc)a)\nn\\
&&+(-1)^{|a||b|}((ba)c-b(ac))\label{twentyone}\\
&&+(a(bc)-(ab)c)\nn\eea
gives the first indication of why $\de^2=0$ is equivalent to the
associativity of $m$. 
\vsp

Next, as in \cite{Ger}, we define a {\bf dot (cup) product}
\[ x\cdot y\]
of cochains $x$, $y\in\CA$ with $D(x)=k$, $D(y)=l$ by
\bea (x\cdot y)(a_1,\dots,a_{k+l})&=&(-1)^{D(x)}\{ m\}\{ x,y\}\{ 
a_1,\dots,a_{k+l}\}\label{twentyfour}\\
&=&(-1)^{kl+|y|(|a_1|+\cdots +|a_k|)}x(a_1,\dots,a_k)
y(a_{k+1},\dots,a_{k+l}),\nn\eea
which is just $m$ on $C^0(A)$ ($k=l=0$). Clearly,
\[ D(x\cdot y)=D(x)+D(y),\;\;\;\; d(x\cdot y)=d(x)+d(y)+1,\;\;\;\;
\mbox{and}\;\;\;\; |x\cdot y|=|x|+|y|.\]

\subsubsection{A new approach: the second level of braces}

\label{seclev}
We can define $\de(x)$ and $x\cdot y$ in terms of the bilinear associative map
$m$ without specifying all the arguments. First, let us take our 
{\bf Z}-graded vector space to be 
\be (B,|\;\; |')=(\CA,D),\label{zgra}\ee
and look at the Hochschild complex $\CB$ where the new adjusted degree of 
homogeneity will be shown by $d'$, the new coupled braces by $\{\; ,\;\}'$,
the new G-bracket by $[\; ,\; ]'$, the new super degree by $|\;\; |'=d+d'$, 
and the new suspended degree by $\|\;\;\|'=|\;\;|'+d'=d$ (mod 2).
This is consistent with~(\ref{zgra}), as
\[ d'(x)=-1,\;\;\; |x|'=d(x)-1=D(x)\;\;\;\mbox{(mod 2),}\]
and
\[ \| x\|'=D(x)-1=d(x)\;\;\;\mbox{(mod 2) for $x\in\CA$.}\]
Note that the original
super degree on $A$ does not have a role at this level. The new
super degree simply comes from the $D$-degree on elements of $B$, similar to
the old super degree on $\CA$ generated by the super degree on $A$. If we 
denote the new suspension operator by $s'$, we will again take
\[ |s'|'=-1\;\;\;\mbox{and}\;\;\; d'(s')=0.\]
We are now ready to define a linear operator $M_1\in C^1(B)$ and 
a bilinear operator $M_2\in C^2(B)$. Let
\be \de(x)=M_1(x)=\{ M_1\}'\{ x\}'\;\defi\; \{ [m,x]\}'\label{d}\ee
(first written in this form by Gerstenhaber in \cite{Ger}), and
\be x\cdot y=M_2(x,y)=\{ M_2\}'\{ x,y\}'\;\defi\; (-1)^{D(x)}\{\{ m\}
\{ x,y\}\}' \label{dot}\ee
(we introduced the second level of braces in Section~\ref{twofour}). Clearly, 
we have
\be |M_1|'=d([m,x])+1-D(x)=d(m)+d(x)+1-d(x)-1=1,\;\;\; 
d'(M_1)=0,\;\;\; \| M_1\|'=1,\ee
and
\be |M_2|'=d(m)+d(x)+d(y)+1-d(x)-1-d(y)-1=0,\;\;\; d'(M_2)=1,\;\;\;
\| M_2\|'=1.\ee
Proposition~\ref{prop2} below works only when ${[m,x]}$ is {\sl
not} modified by a sign depending on $x$, hence we adopt~(\ref{d}) as the
definition of the Hochschild differential instead of the 
classical~(\ref{eighteen}). 
We will revisit these ideas by constructing a strongly homotopy associative
product~$M$ on $\CB$ as in \cite{Get1}, starting from a strongly homotopy 
associative structure $m\in\CA$. 

\subsubsection{Properties of the differential and the dot product}
\label{AA}

We can summarize several properties of the differential and the dot
product as follows:

\begin{thm} $(\CA,M_1,M_2)$ is a differential graded associative algebra.
\end{thm}

The Theorem first of all asserts that $M_2$ is associative, or
\[ \{ M_2\}'\{ M_2\}'=M_2\circ M_2=0.\]
Note that the factor $(-1)^{D(x)}$ is absolutely
necessary for the associativity of the dot product:

\begin{prop} For the dot product~(\ref{dot}), we have
\[ (x\cdot y)\cdot z-x\cdot(y\cdot z)=(-1)^{D(y)}
\{\{ m\circ m\}\{ x,y,z\}\}' ;\]
associativity of $M_2$ follows from that of $m$.
\end{prop}

{\em Proof.} The left hand side is 
\bea &&(-1)^{D(x\cdot y)}\{ m\}\{ x\cdot y,z\}-(-1)^{D(x)}\{ m\}\{ x,y\cdot 
z\}\nn\\
&=&(-1)^{D(x)+D(y)+D(x)}\{ m\}\{\{ m\}\{ x,y\} ,z\} -(-1)^{D(x)+D(y)}\{ m\}
\{ x,\{ m\}\{ y,z\}\}\nn\\
&=&(-1)^{D(y)}(\{ m\}\{\{ m\}\{ x,y\} ,z\} +(-1)^{d(x)}\{ m\}\{ x,\{ m\}
\{ y,z\}\})\nn\\
&=&(-1)^{D(y)}\{ m\}\{ m\}\{ x,y,z\} ,\nn\eea
which is equal to the right hand side. This proof is a generalization
of~(\ref{twentythree})! $\Box$
\vsp

Secondly, we understand that 
\[ \de^2=M_1^2=\frac{1}{2}[M_1,M_1]'=  0.\]

\begin{prop}\label{dsq}
We have $\de^2=0$ if and only if $m\circ m=0$.
\end{prop}

{\em Proof.} We take $|m|=0$ and $d(m)=1$, and compute
\bea {[m,[m,x]]}=&&[[m,m],x]+(-1)^{d(m)d(m)+|m||m|}[m,[m,x]]\nn\\
=&&2[m\circ m,x]-[m,[m,x]],\nn\eea
which means
\[ \de^2(x)=\{ [m\circ m,x]\}'\;\;\;\;\forall x.\]
Clearly associativity implies $\de^2=0$. Conversely, by setting $\de^2=0$
and $x=\mbox{id}$, we obtain
\[ {[m\circ m,\mbox{id}]}\{ a,b,c\} =2((ab)c-a(bc))=0,\]
or $m\circ m=0$. 
(Penkava has a similar proof in \cite{Pen1}.)$\Box$
\vsp

Gerstenhaber shows in \cite{Ger} that $\de$ is a derivation of 
the dot product (with respect to the $D$-grading on $B=\CA$).
In other words, we have
\be M_1\circ M_2-M_2\circ M_1={[M_1,M_2]'}=0.\ee
We can furthermore prove

\begin{prop} \label{prop2} For any algebra $A$ with a bilinear product $m$ 
and a dot product and differential defined as above in terms of $m$, we have
\[\de(x\cdot y)-\de(x)\cdot y-(-1)^{D(x)}x\cdot\de(y)
=(-1)^{D(x)}\{\{ m\circ m\}\{ x,y\}\}' .\]
In particular, $\de =M_1$ is a derivation of $M_2$ if and only if $m$ is
associative.
\end{prop}

{\em Proof.} We write the left hand side as
\bea \label{bal}&&(-1)^{D(x)}[m,\{ m\}\{ x,y\} ]\\
&&-(-1)^{D(x)+1}\{ m\}\{ [m,x],y\} -(-1)^{D(x)+D(x)}\{ x,[m,y]\}\nn\\
&=&(-1)^{D(x)}(\{ m\}\{\{ m\}\{ x,y\}\} +(-1)^{d(x)+d(y)}
\{ m\}\{ x,y\}\{ m\}\nn\\
&&+\{ m\}\{ [m,x],y\} +(-1)^{d(x)}\{ m\}\{ x,[m,y]\} ).\nn\eea
Meanwhile, substituting $m$, $m$, $x$, $y$, for $x$, $y$, $z_1$,
$z_2$ respectively in Lemma~\ref{triple}, we obtain
\bea &&\{ m\}\{ m\}\{ x,y\}\nn\\
&=&\{ m\}\{\{ m\}\{ x,y\}\} +\{ m\}\{ [m,x],y\}\nn\\
&&+(-1)^{d(x)}\{ m\}\{ x,[m,y]\} +(-1)^{d(x)+d(y)}\{ m\}\{ x,y\}\{ m\} .
\nn\eea
Then the right hand side of~(\ref{bal}) must be equal to
\[ (-1)^{D(x)}\{\{ m\circ m\}\{ x,y\}\}' .\]
$\Box$
\vsp
 
\subsubsection{Bialgebra cohomology}

The differential and the dot product can be defined on $\CBB$
for an associative algebra $A$
via exactly the same formulas as above, thanks to the existence of
composition and Gerstenhaber bracket on the extended complex. 
Still, the extended complex, or more precisely its subcomplex
\be \tC ,\label{subc} \ee
is more useful in the context of bialgebras (see Gerstenhaber and Schack 
\cite{GeS} and Stasheff \cite{Jim}). We first recall that the cohomology 
complex for a coassociative coalgebra $(A,\com)$ (with comodule $A$) is
\be \bC, \ee
and the differential $\bde :\bar{C}^k\ra\bar{C}^{k+1}$ is given by
\be \bde(x)=[\com,x]\ee
(our interpretation). Since we recognize the condition for coassociativity as
\be \com\circ\com =0,\label{coassoc}\ee
we have 
\be \bde^2=0.\ee
Again we have to reconcile (\ref{coassoc}) with the
well-known condition
\[ (\com\ot\ide -\ide\ot\com)\circ\com =0,\]
but then $\com$ applied to a tensor $a\ot b$ in the image of $\com$ is
by definition
\bea \{\com\}\{ a,b\} &=&\{\{\com\}\{ a\},b\} +(-1)^{d(\com)d(a)+|\com |
|a|}\{ a,\{\com\}\{ b\}\}\nn\\
&=&\{\{\com\}\{ a\},b\} -\{ a,\{\com\}\{ b\}\}\nn\\
&=&(\com\ot\ide -\ide\ot\com)(a,b).\nn\eea
Note that we have for $x:A\ra A^{\ot j}$ ($j\geq 1$)
\bea &&D(\{\com\}\{ x\})=D(x)+\ma\{ 0,D(\com)-R(x)\} =1+\ma\{ 0,1-j\}
=1\nn\\
&&D(\{ x\}\{\com\})=D(\com)+\ma\{ 0,D(x)-R(\com)\} =1+\ma\{ 0,1-2\} =1\nn\\
&&R(\{\com\}\{ x\})=R(\com)+\ma\{ 0,R(x)-D(\com)\} =2+\ma\{ 0,j-1\} =j+1\nn\\
&&R(\{ x\}\{\com\})=R(x)+\ma\{ 0,R(\com)-D(x)\} =j+\ma\{ 0,2-1\}
=j+1,\nn\eea
so that 
\be D([\com,x])=D(x)=1\;\;\;{\rm and}\;\;\;R([\com,x])=R(x)+1=j+1;\ee
the differential increases the degree of a cochain by one.
\vsp

In \cite{GeS} Gerstenhaber and Schack define the cohomology differential $\hde$
for a bialgebra $A$ on~(\ref{subc}) (in fact, even more generally for any
birepresentation of this bialgebra) as a signed sum of algebra and
coalgebra differentials. First, we need a grading on $\tC$, namely
\[ \hat{C}^{\bullet}(A)=\op_{n\geq -1}\hat{C}^n(A),\]
with
\be \hat{C}^n=\op_{i+j=n+1}Hom(A^{\ot i};A^{\ot j}).\ee
Then if $x\in Hom(A^{\ot i};A^{\ot j})$, we verify that neither ${[m,x]}$
nor ${[\com,x]}$ stays completely in $\hat{C}^{n+1}(A)$ by computing
\bea &&D(\{ m\}\{ x\})=i+2\de_{j,0}+\de_{j,1}  \nn\\
&&R(\{ m\}\{ x\})=j-1+2\de_{j,0}+\de_{j,1}  \nn\\
&&D(\{ x\}\{ m\})=i+1+\de_{i,0}   \nn\\
&&R(\{ x\}\{ m\})=j+\de_{i,0}   \nn\eea
and
\bea &&D(\{\com\}\{ x\})=i+\de_{j,0}   \nn\\
&&R(\{\com\}\{ x\})=j+1+\de_{j,0}   \nn\\
&&D(\{ x\}\{\com\})=i-1+2\de_{i,0}+\de_{i,1}    \nn\\
&&R(\{ x\}\{\com\})=j+2\de_{i,0}+\de_{i,1}     .\nn\eea
The differential $\hde$ is defined naturally on $Hom(A^{\ot i};A^{\ot j})$
as the signed sum of the algebra cohomology differential (for the $A$-module
$A^{\ot j}$) and the coalgebra cohomology differential (for the $A$-comodule
$A^{\ot i}$). We decode this statement as follows (see Giaquinto's thesis
\cite{Gia} for very clear definitions).
\vsp

For a bialgebra $A$, we can define a left $A$-module structure 
\[ m_L:A\ot A^{\ot n}\ra A^{\ot n}\]
and a left $A$-comodule structure
\[ \com_L:A^{\ot n}\ra A\ot A^{\ot n}\]
on $A^{\ot n}$. When $n=1$, $m_L$ and $\com_L$ are just the multiplication
and comultiplication maps. When $n=2$, we have
\bea m_L(a,b\ot c)&=&\com(a)\cdot(b\ot c)\nn\\
&=&\sum_{\alpha}\{ a_{\alpha}^{(1)},a_{\alpha}^{(2)}\}\cdot\{ b,c\}\nn\\
&=&\pm\sum_{\alpha}\{ a_{\alpha}^{(1)}b,a_{\alpha}^{(2)}c\}\nn\\
&=&\{ m,m\}\{\sigma_2\}\{\com,\ide,\ide\}\{ a,b,c\}\nn\\
&=&\{ m^{\ot 2}\}\{\sigma_2\}\{\com,\ide^{\ot 2}\}\{ a,b,c\} ,\nn\eea
where $\sigma_2:A^{\ot 4}\ra A^{\ot 4}$ is the signed permutation of tensor
factors given by
\[ \sigma_2=\left( \begin{array}{cccc} 1 & 2 & 3 & 4 \\ 1 & 3 & 2 & 4
\end{array}\right) .\]
Similarly, 
\bea \com_L(a\ot b)&=&\pm\{ m,\ide,\ide\}\{\tau_2\}\{\com(a),\com(b)\}\nn\\
&=&\{ m,\ide^{\ot 2}\}\{\tau_2\}\{\com^{\ot 2}\}\{ a,b\} ,\nn\eea
or
\[ \com_L(a\ot b)=\pm\sum_{\alpha,\beta}\{ a_{\alpha}^{(1)}b_{\beta}^{(1)},
a_{\alpha}^{(2)},b_{\beta}^{(2)}\} ,\]
where $\tau_2=\sigma_2^{-1}=\sigma_2$. We proceed by induction, and obtain
\bea &&m_L(a,b_1\ot\cdots\ot b_n)\\
&=&\{ m^{\ot n}\}\{\sigma_n\}\{\com,\ide^{\ot(2n-2)}\}\{\com,\ide^{\ot(2n-3)}
\}\cdots\{\com,\ide^{\ot n}\}\{ a,b_1,\dots,a_n\} \nn\eea
and
\bea &&\com_L(a_1\ot\cdots\ot a_n)\\
&=&\{ m,\ide^{\ot n}\}\{ m,\ide^{\ot(n+1)}\}\cdots\{ m,\ide^{\ot(2n-2)}\}
\{\tau_n\}\{\com^{\ot n}\}\{ a_1,\dots,a_n\} ,\nn\eea
where $\sigma_n,\tau_n:A^{\ot 2n}\ra A^{\ot 2n}$ are the signed
permutations given by
\[ \sigma_n=\left( \begin{array}{ccccccc} 1 & 2 & 3 & 4 & \dots & 2n-1 &
2n \\ 1 & n+1 & 2 & n+2 & \dots & n & 2n \end{array}\right) \]
and
\[ \tau_n=\sigma_n^{-1}=\left( \begin{array}{cccccccccc} 1 & 2 & 3 & \dots
& n & n+1 & n+2 & n+3 & \dots & 2n \\ 1 & 3 & 5 & \dots & 2n-1 & 2 & 4 & 6 &
\dots & 2n \end{array}\right) .\]
On the other hand, the right module and comodule structure maps on $A^{\ot n}$
are given by
\bea &&m_R(b_1\ot\cdots\ot b_n,a)\\
&=&\{ m^{\ot n}\}\{\sigma_n\}\{\ide^{\ot(2n-2)},\com\}\cdots\{\ide^{\ot n},
\com\}\{ b_1,\dots,b_n,a\}\nn\eea
and
\bea &&\com_R(a_1\ot\cdots\ot a_n)\\
&=&\{\ide^{\ot n},m\}\cdots\{\ide^{\ot(2n-2)},m\}\{\tau_n\}\{\com^{\ot n}\}
\{ a_1,\dots,a_n\} .\nn\eea
Then the algebra cohomology differential
\[ \de :Hom(A^{\ot i};A^{\ot j})\ra Hom(A^{\ot(i+1)};A^{\ot j})\]
for the $A$-bimodule $A^{\ot j}$ can be written as
\be \de(x)=\pm\{ m_L\}\{\ide,x\}\pm\{ x\circ m\} +\{ m_R\}\{ x,\ide\} ,\ee
and the coalgebra cohomology differential
\[ \bde :Hom(A^{\ot i};A^{\ot j})\ra Hom(A^{\ot i};A^{\ot(j+1)})\]
for the $A$-bicomodule $A^{\ot i}$ is
\be \bde(x)=\pm\{\ide,x\}\{\com_L\} +\{\com\circ x\}\pm\{
x,\ide\}\{\com_R\} .\ee
The definition of the bialgebra cohomology differential
\[ \hde :\hat{C}^n(A)\ra\hat{C}^{n+1}(A) \]
on $Hom(A^{\ot i};A^{\ot j})\subset\hat{C}^n(A)$ is then
\be \hde(x)=\de(x) +\bde(x),\ee
and it is known to be square-zero as the two differentials commute.

\subsection{Gerstenhaber and Batalin-Vilkovisky algebras}
\label{CC}

In \cite{VG}, three groups of identities on $\CA$ satisfied by the braces 
(\ref{four}),
the dot product $M_2$, and the differential $\de=M_1$ are 
singled out as the
definition of a {\bf homotopy G-algebra} ($G$ for Gerstenhaber). These are:
(i) the higher pre-Jacobi identities (\ref{eight}), (ii) the distributivity
of $M$ over the braces, namely
\be \{ x_1\cdot x_2\}\{ y_1,\dots,y_n\} =\sum_{k=0}^{n}
(-1)^{D(x_2)(d(y_1)+\cdots +d(y_k))+|x_2|(|y_1|+\cdots +|y_k|)}
\{ x_1\}\{ y_1,\dots,
y_k\}\cdot\{ x_2\}\{ y_{k+1},\dots,y_n\} ,\label{twentyseven}\ee
and (iii)
\bea &&\de(\{ x\}\{ y_1,\dots,y_{n+1}\} )-\{ \de(x)\}\{ y_1,\dots,y_{n+1}\}
\label{twentyeight}\\
&&-(-1)^{d(x)}\sum_{i=1}^{n+1}(-1)^{d(y_1)+\cdots +d(y_{i-1})}\{ x\}
\{ y_1,\dots,\de(y_i),\dots,y_{n+1}\}\nn\\
&=&(-1)^{D(x)d(y_1)+|y_1||x|}y_1\cdot\{ x\}\{ y_2,\dots,y_{n+1}\}\nn\\
&&-(-1)^{d(x)}\sum_{i=1}^n(-1)^{d(y_1)+\dots +d(y_i)}\{ x\}\{ y_1,\dots,
y_i\cdot y_{i+1},\dots,y_{n+1}\}\nn\\
&&+(-1)^{d(x)+d(y_1)+\cdots +d(y_n)}\{ x\}\{ y_1,\dots,y_n\}\cdot y_{n+1},
\nn\eea
a higher homotopy identity. We see~(\ref{twentyseven}) as a special case of 
the higher pre-Jacobi identity:
\bea &&\{ x_1\cdot x_2\}\{ y_1,\dots,y_n\}\nn\\ 
&=&(-1)^{D(x_1)}\{ m\}\{ x_1,x_2\}\{ y_1,\dots,y_n\}\nn\\
&=&(-1)^{D(x_1)}\sum_{k=0}^n(-1)^{d(x_2)(d(y_1)+\cdots +d(y_k))+|x_2|(|y_1|+
\cdots +|y_k|)}\mxy\nn\\
&=&\sum_{k=0}^n(-1)^{d(x_1)+1+D(x_2)(d(y_1)+\cdots +d(y_k))+d(y_1)+\cdots
+d(y_k)+{\rm super}}\mxy\nn\\
&=&\sum_{k=0}^n(-1)^{D(\{ x_1\}\{ y_1,\dots,y_k\})+D(x_2)(d(y_1)+\cdots +
d(y_k))+{\rm super}}\mxy\nn\\
&=&\sum_{k=0}^n(-1)^{D(x_2)(d(y_1)+\cdots +d(y_k))+{\rm super}}
\{ x_1\}\{ y_1,\dots,y_k\}\cdot\{ x_2\}\{ y_{k+1},\dots,y_n\} .\nn\eea
The identity (\ref{twentyeight}) can again be unraveled by writing the terms in
\[ \de(\{ x\}\{ y_1,\dots,y_{n+1}\})=[m,\{ x\}\{ y_1,\dots,y_{n+1}\} ]\]
explicitly.
\vsp

An ordinary {\bf G-algebra} is, on the other hand, a graded commutative and 
associative algebra $(A,|\;\; |)$
(with a bilinear map called a ``dot product'') together with an odd ``Poisson
bracket'' $\{\; ,\;\}$ satisfying identities similar to those in $\CA$. 
Namely, we have

(i) {\em antisymmetry in the associated, suspended-graded Lie algebra}
$\hat{A}\stackrel{\rm def}{=}\sum_jA^{j-1}$:
\[ \{ a,b\} =-(-1)^{(|a|-1)(|b|-1)}\{ b,a\} ;\]

(ii) {\em the Leibniz rule, or suspended-graded derivation property,
in} $\hat{A}$:
\[ \{ a,\{ b,c\}\} =\{\{ a,b\} ,c\} +(-1)^{(|a|-1)(|b|-1)}\{ b,\{ a,c\}\} ;\]
and

(iii) {\em the graded derivation (Poisson) rule with respect to the dot product
in A}:
\[ \{ a,b\cdot c\} =\{ a,b\}\cdot c+(-1)^{(|a|-1)|b|}b\cdot\{ a,c\} .\]

Proposition~\ref{gers} and Lemmas \ref{A} and \ref{B} below give us the
prime example of a G-algebra:

\begin{prop}[Gerstenhaber] The cohomology $H(\CA,\de)$ of the 
Hochschild complex with
the induced dot product $M_2$ and Gerstenhaber bracket $[\;,\; ]$
has the structure of a G-algebra.\label{gers}
\end{prop}

{\it Proof.} We partially follow suggestions in \cite{VG}.
For simplicity, we assume $A$ has no original super grading 
(otherwise we will have to modify all statements according to the bigrading
$(D,|\;\; |)$). First of all, $(\CA,M_2)$ has been shown to be a $D$-graded
associative algebra, and $[\;,\; ]$ is a $D$-odd bracket:
\[ D([x,y])-D(x)-D(y)=d(x)+d(y)+1-d(x)-1-d(y)-1=-1.\]
Secondly, $M_2$ commutes with $M_1=\de$, and hence descends to the 
$\de$-cohomology. Furthermore, $M_2$ is homotopy commutative:
\[ x\cdot y-(-1)^{D(x)D(y)}y\cdot x=(-1)^{d(x)}(\de(x\circ y)-\de(x)\circ y
-(-1)^{d(x)}x\circ\de(y))\] 
from (\ref{twentyeight}) with $n=0$.
We show that $\de$ is also a derivation of the
$G$-bracket with respect to the $d=D-1$-grading, i.e.
\[ \de([x,y])-[\de(x),y]-(-1)^{d(x)}[x,\de(y)]=0,\]
in Lemma~\ref{A} below (thus the bracket is defined on the cohomology).
We have seen that the $G$-bracket satisfies the
graded antisymmetry and graded Leibniz (derivation) properties with respect to
$d$. Finally, we show the homotopy Poisson rule in Lemma~\ref{B}.$\Box$

\begin{lemma}\label{A} $\de([x,y])-[\de(x),y]-(-1)^{d(x)}[x,\de(y)]=0.$
\end{lemma}

{\it Proof.} This is a direct result of the definition of $\de$ and
the derivation property of the $G$-bracket: the left hand side is exactly
\bea && {[m,[x,y]\; ]}-([\; [m,x],y]+(-1)^{d(x)}[x,[m,y]\; ])\nn\\
&=&{[m,[x,y]\; ]}-([\; [m,x],y]+(-1)^{d(x)d(m)+|x||m|}[x,[m,y]\; ])\nn\\
&=&0.\nn\eea
Note that the result is true even when $A$ does have a super grading.$\Box$

\begin{lemma}\label{B} Ignoring the super grading, we have
\bea &&{[x,y\cdot z]}-[x,y]\cdot z-(-1)^{d(x)D(y)}y\cdot [x,z]\nn\\
&=&(-1)^{d(x)+D(y)}(\de(\{ x\}\{ y,z\})-\{\de(x)\}\{ y,z\} -(-1)^{d(x)}
\{ x\}\{\de(y),z\} -(-1)^{d(x)+d(y)}\{ x\}\{ y,\de(z)\}).\nn\eea
\end{lemma}

{\it Proof.} Appendix.$\Box$
\vsp

Many examples of $G$-algebras are in fact {\bf Batalin-Vilkovisky (BV) 
algebras} \cite{BV},
where the bracket $\{ \; ,\;\}$ is obtained from an odd, square zero,
second order differential operator $\tri$ on a supercommutative and associative
algebra $A$ (we steer away from duplicate notation by omitting the dot and
replacing the curly braces by $\{\; ,\;\}_{\tri}$ from here on): the braces
\bea \{ a,b\}_{\tri} &\stackrel{\rm def}{=}&(-1)^{|a|}\Phi_{\tri}^2(a,b)\nn\\
&=&(-1)^{|a|}\tri(ab)-(-1)^{|a|}\tri(a)b-a\tri(b)\label{bvb}\eea
measure the deviation of $\tri$ from being a first order differential
operator. In \cite{A} the notion
was generalized to an arbitrary algebra and an arbitrary linear operator,
and analogues of the above properties of the dot product and the bracket were 
discussed. With our new language, we can write the (generalized) BV bracket 
as
\be\label{bv2}\{ a,b\}_{\tri}=(-1)^{|a|-1}\{ s\}[m,\tri ]\{ a,b\} =(-1)^{|a|-1}
\de(\tri)(a,b) ,\ee
and prove its properties in a few lines, in a major change from the usual
methods (compare with \cite{A} and the references therein). We have
\[ d(\{\; ,\;\}_{\tri})=d(\tri)+d(m)=1\;\;\;\mbox{and}\;\;\; |\{\; ,\;\}_{\tri}|
=|\tri |+|m|={\rm odd},\]
so that
\[ \| \{\; ,\;\}_{\tri}\| ={\rm even}.\]
Experience shows that it is better to treat $[m,\tri ]$ as a bilinear operator
on $A$, and $\{\; ,\;\}_{\tri}$ as a bilinear operator
on the suspended-graded space $sA$. Recall that the Gerstenhaber 
bracket between suspended-graded operators is
suspended-graded antisymmetric. We now rewrite the following statement
in our new notation and give a proof which uses nothing deeper than the
definition of the $\Phi$'s and Lemma~\ref{yasa} in Section~\ref{phi}.
\vsp

\begin{prop}[Akman \cite{A}] For a superalgebra $(A,m)$ and an odd linear 
operator $\tri$ on $A$, the BV bracket defined by~(\ref{bvb}) satisfies the 
following properties.

(i) Modified $\|\;\|$-graded antisymmetry:
\[ \{ sa,sb\}_{\tri}+(-1)^{\| a\|\,\| b\|}\{ sb,sa\}_{\tri} 
=(-1)^{|a|}\{ s\}
\{\tilde{\Phi}_{\tri}^2\}\{ a,b\} =\{ sa,sb\}_{\tri}^{\tilde{}}.\]
(ii) Modified $\|\;\;\|$-graded Leibniz rule:
\bea \{ sa,\{ sb,sc\}_{\tri}\}_{\tri}&-&
\{\{ sa,sb\}_{\tri},sc\}_{\tri}-(-1)^{\| a\|\,\| b\|}\{ sb,\{ sa,sc\}_{\tri}
\}_{\tri}
\nn\\ &=& (-1)^{|b|}\{ s\}\{\Phi_{\tri^2}^3-[\tri,\Phi_{\tri}^3]\}\{ a,b,c\} .
\nn\eea
(iii) Modified Poisson rule:
\[ \{ sa,s(bc)\}_{\tri}-\{ sa,sb\}_{\tri}c-(-1)^{\| a\|\, |b|}b
\{ sa,sc\}_{\tri} 
=(-1)^{|a|}\{ s\}\{\Phi_{\tri}^3\}
\{ a,b,c\} .\]
(iv) Modified derivation rule for $\tri$:
\[ \{\tri\}\{ sa,sb\}_{\tri}-\{\{\tri\}\{ sa\} ,sb\}_{\tri}
-(-1)^{\| a\|}\{ sa,\{\tri\}\{ sb\}\}_{\tri} 
=(-1)^{\| a\|}\{ s\}\{\Phi_{\tri^2}^2\}\{ a,b\}=\{ sa,sb\}_{\tri^2} .\]
The tilde in Property~(i) says that the 
bilinear product $m$ is replaced by its super 
antisymmetrization~(\ref{sixteen}), which we will denote by $l$.
\end{prop}

{\it Proof.} (i) We have
\bea &&\{ s^{-1}\}(\{ sa,sb\}_{\tri}+(-1)^{\| a\|\,\| b\|}\{
sb,sa\}_{\tri})\nn\\
&=&(-1)^{\| a\|}\mD\{ a,b\} 
+(-1)^{\| a\|\,\| b\| +\| b\|}\mD\{ b,a\}\nn\\
&=&(-1)^{\| a\|}(\mD\{ a,b\} -(-1)^{|a||b|}\mD\{ b,a\} )\nn\\
&=&(-1)^{\| a\|}\mD\{ a\}\{ b\}\nn\\
&=&(-1)^{\| a\|}[l,\tri ]\{ a,b\} \nn\\
&=&(-1)^{|a|}\{\tilde{\Phi}_{\tri}^2\}\{ a,b\}\nn\\ &=&\{ s^{-1}\}
\{ sa,sb\}_{\tri}^{\tilde{}}.\nn\eea

(ii) This is in fact the third identity in Lemma~\ref{yasa}. First, we have
\bea &&\{ s^{-1}\}(\{ sa,\{ sb,sc\}_{\tri}\}_{\tri}-
\{\{ sa,sb\}_{\tri},sc\}_{\tri}
-(-1)^{\| a\|\,\| b\|}\{ sb,\{ sa,sc\}_{\tri}\}_{\tri})\nn\\
&=&(-1)^{\| a\| +\| b\|}\Phi_{\tri}^2(a,\Phi_{\tri}^2(b,c))-(-1)^{\| a\| +
\| a\| +\| b\|}\Phi_{\tri}^2(\Phi_{\tri}^2(a,b),c)-(-1)^{\| a\|\,\| b\| +
\| b\| +\| a\|}\Phi_{\tri}^2(b,\Phi_{\tri}^2(a,c))\nn\\
&=&(-1)^{|b|}[\Phi_{\tri}^2,\adj(\Phi_{\tri}^2)\{ a\} ](b,c)\nn\eea
by definition of BV and G-brackets. But by the Lemma (where $T=U=\tri$)
this is exactly
\[ (-1)^{|b|}(\Phi_{\tri^2}^3(a,b,c)-[\tri,\Phi_{\tri}^3](a,b,c));\]
note that $\tri^2=\frac{1}{2}[\tri,\tri ]$.

(iii) This is just the definition of $\Phi_{\tri}^3$:
\bea &&(-1)^{|a|}\{\Phi_{\tri}^3\}\{ a,b,c\} \nn\\
&=&(-1)^{|a|}(\{\Phi_{\tri}^2\}\{ a,bc\} 
-\{\Phi_{\tri}^2\}\{ a,b\} c-(-1)^{|b|(|a|-1)}
b\{\Phi_{\tri}^2\}\{ a,c\} )\nn\\
&=&\{ s^{-1}\}
(\{ sa,s(bc)\}_{\tri}-\{ sa,sb\}_{\tri}c-(-1)^{\| a\| |b|}b\{ sa,sc\}_{\tri}).
\nn\eea

(iv) The left hand side is given by
\bea &&\{ s^{-1}\}{\rm LHS}\nn\\
&=&(-1)^{\| a\|}\tri(\Phi_{\tri}^2(a,b)-(-1)^{\| a\| +1}\Phi_{\tri}^2(\tri(a),
b)-(-1)^{\| a\| +\| a\|}\Phi_{\tri}^2(a,\tri(b))\nn\\
&=&(-1)^{\| a\|}([\tri,\Phi_{\tri}^2](a,b))\nn\\
&=&(-1)^{\| a\|}\Phi_{\tri^2}^2(a,b)\nn\eea
by Lemma~\ref{yasa}.$\Box$

\subsection{Strongly homotopy associative algebras}

\subsubsection{Definition}
\label{fare}

We are now back to $\CA =Hom(TA;A)$.
Strongly homotopy associative ($\Ai$) algebras were introduced by Stasheff
in \cite{St}. We will partially follow Getzler's approach in \cite{Get1}.
See also \cite{Pen1}, \cite{Pen2}, and \cite{PS}.
The associative bilinear product $m\in C^2(A)$ in Section~\ref{onetwo}
is now replaced with the formal sum
\be m=m_1+m_2+\cdots\label{twentynine}\ee
of multilinear products
\[ m_k:A^{\ot k}\ra A.\]
Like Getzler, we will define an {\bf $\Ai$ 
algebra} to be a super
({\bf Z}) graded vector space $A$ with some cochain $m\in\CA$ satisfying
\be \tilde{m}\circ\tilde{m}=0\label{thirty}\ee
in addition to the parity conditions
\be (-1)^{|m_k|}=(-1)^k,\;\;\;\; k\geq 1,\label{thirtyone}\ee
so that
\be (-1)^{\| m_k\|}=-1,\ee
and
\be {[\tm,\tm ]}=2\,\tm\circ\tm =0.\ee

The condition (\ref{thirty}) 
makes $T(sA)$ into a {\sl differential} graded coalgebra with respect to
the suspended grading;
we look for (homotopy) associativity and other desirable properties in the 
unadjusted products $m_n$ on $TA$ with the bigrading. This master identity 
unfolds as
\be \sum_{i+j=n+1}\tilde{m}_i\circ\tilde{m}_j=0\;\;\;\mbox{for each 
$n\geq 1$,}\ee
or as
\be \sum_{i+j=n+1}[\tm_i,\tm_j]=0\;\;\;\mbox{for each $n\geq 1$}.\ee
Equivalently, we may write
\be \sum_{i+j=n+1}\sum_{k=0}^{i-1}(-1)^{\| m_j\| (\| a_1\| +\cdots
+\| a_k\| )}
\;\tilde{m}_i(a_1,\dots,a_k,\tilde{m}_j(a_{k+1},\dots),\dots,a_n)=0
\;\;\;\;\mbox{for all $n$.}\label{thirtytwo}\ee

\begin{prop}
The statement $\tilde{m}\circ\tilde{m}=0$ is equivalent to
\be\label{vay} \sum_{i+j=n+1}\sum_{k=0}^{i-1}(-1)^{j(|a_1|+\cdots+|a_k|)
+jk+j+k}\; m_i(a_1,\dots,a_k,m_j(a_{k+1},\dots,a_{k+j}),a_{k+j+1},\dots,a_n)=0
\ee
for all $n\geq 1$, similar to the original $\Ai$ identity in \cite{LS}.
\end{prop}

\begin{Rm} See Markl's explanation of sign discrepancy in \cite{Ma}, 
Example~1.6.
\end{Rm}

{\it Proof.} From (\ref{thirtytwo}) and the definition of $\tilde{m}$, 
the power of $(-1)$ in front of
\[ m_i(a_1,\dots,a_k,m_j(a_{k+1},\dots,a_{k+j}),a_{k+j+1},\dots,a_n)\]
is
\bea &&\sum_{t=1}^k\| a_t\| +\sum_{t=1}^j(j-t)\| a_{k+t}\| +\sum_{t=1}^k
(i-t)\| a_t\|\nn\\
&&+(i-k-1)(1+\sum_{t=1}^j\| a_{k+t}\| )+\sum_{t=1}^{n-k-j}(i-k-1-t)\|
a_{k+j+t}\|\nn\\
&=&\sum_{t=1}^k(i+j-1-t)\| a_t\| +j\sum_{t+1}^k\| a_t\| +\sum_{t=1}^j
(i+j-1-t-k)\| a_{k+t}\|\nn\\
&&+\sum_{t=1}^{n-k-j}(i-k-1-t)\| a_{k+j+t}\| +i+k+1\;\;\;\mbox{(mod 2)}\nn\\
&=&\sum_{t=1}^n(n-t)\| a_t\| +j\sum_{t=1}^k\| a_t\| +i+k+1\;\;\;
\mbox{(mod 2)}\nn\\
&=&\sum_{t=1}^n(n-t)\| a_t\| +j\sum_{t=1}^k|a_t|+jk+i+k+1\;\;\;
\mbox{(mod 2)}\nn\\
&=&\sum_{t=1}^n(n-t)\| a_t\| +j\sum_{t=1}^k|a_t|+jk+j+n+k\;\;\;
\mbox{(mod 2)}.\nn\eea
Since the first term and $n$ are independent of 
$i$, $j$, and $k$, we are done.$\Box$ 
\vsp

\begin{prop} \label{eqq}Another equivalent statement is
\be \label{mbm} \sum_{i+j=n+1}(-1)^i[m_i,m_j]=2\sum_{i+j=n+1}(-1)^im_i\circ
m_j=0,\ee
or
\[ \sum_{i+j=n+1}(-1)^j[m_i,m_j]=2\sum_{i+j=n+1}(-1)^jm_i\circ m_j=0.\]
\end{prop}

{\it Proof.} Appendix.$\Box$
\vsp

The lowest identities in (\ref{mbm}) have been studied
extensively. We have
\be m_1^2(a)=0\;\;\;\;\forall a\in A\label{thirtythree}\ee
for $n=1$, which says $m_1$ is an odd differential on $A$. For $n=2$,
\be m_1(m_2(a,b))-m_2(m_1(a),b)-(-1)^{|a|}m_2(a,m_1(b))=0\label{thirtyfour}\ee
says that $m_1$ is a derivation of the even bilinear product $m_2$. We have
\bea &&m_1(m_3(a,b,c))\nn\\
&&-(m_2(m_2(a,b),c)-m_2(a,m_2(b,c)))\label{thirtyfive}\\
&&+m_3(m_1(a),b,c)+(-1)^{|a|}m_3(a,m_1(b),c)+(-1)^{|a|+|b|}m_3(a,b,m_1(c))\nn\\
&&=0\nn\eea
when $n=3$, indicating that $m_2$ is ``associative up to homotopy''
(when $m=m_2$ only, $\tilde{m}\circ\tilde{m}=0$ is equivalent to $m\circ m=0$).
If all but finitely many $m_k$ are identically zero, we obtain some 
well-known structures,
such as a differential complex for $m_1\neq 0$, an associative algebra for
$m_2\neq 0$, and a differential graded associative algebra for $m_1$, 
$m_2\neq 0$.
\vsp

\begin{ex} The maps $M_1$ and $M_2$ on $\CA$ define a differential graded
algebra (or a truncated $\Ai$ algebra) structure.
\end{ex}

\begin{ex} A nontrivial $\Ai$ algebra structure can be defined on any 
associative 
algebra by setting $m_n=0$ for odd $n$ and the (unambiguous) $n$-fold product
for even $n$ (especially in an algebra with no super grading, we do not
expect to have nonzero odd multilinear operators). 
In~(\ref{mbm}) the only nonzero expressions will be
for $n=3$, 5, 7, ... ($n+1=4$, 6, 8, ...) where we have $\tilde{m}_i\circ
\tilde{m}_j$ terms only for $i$, $j$ both even, adding up to $n+1$. 
Then fixing $i$, $j$ as above, we get the alternating expression
$(-1)^k$ as the coefficient of 
\[ m_i(a_1,\dots,a_k,m_j(a_{k+1},\dots,a_{k+j}),a_{k+j+1},\dots,a_n)
=a_1\cdots a_n,\]
and there are an even number of $k$'s. Note how this approach differs from
the one in the previous example: the original super degree on $A$ is a 
hindrance.
\end{ex}

\subsubsection{Properties of the differential and the higher products}
\label{BB}

We go on to define an $\Ai$ product $M\in\CB$, with $\de =M_1$,
a Hochschild-type differential. It seems that this time
\be \de(\tx)=[\tm,\tx ]\label{thirtysix}\ee
is a good candidate, because the proof of the statement
\be \de^2=0\label{thirtyseven}\ee
follows that of Proposition~\ref{dsq}: although $\tilde{m}$ is now a formal
infinite sum, the expression $[\tilde{m},\tilde{m}]$ is again $2\,\tilde{m}
\circ\tilde{m}$. The grading $\|\tilde{m}_k\|$ which is uniformly odd for
all $k$ makes it possible for us to treat $\tilde{m}$  as one quantity when
it comes to writing out the Gerstenhaber bracket, in contrast to $m$ and
its varying double grading.

\begin{Rm} We may think of an $\Ai$~algebra as possessing a product 
$\tilde{m}$ which is a derivation of itself, in the sense of 
Section~\ref{phi}. Generalizations of $\Ai$~algebras can then be obtained
by producing an $\tilde{m}$ which is a higher order differential operator
with respect to itself!
\end{Rm}

The operators $M_1$ and $M_2$ on the old complex (\ref{two}) is generalized to
\[ M=M_1+M_2+\cdots\in\CB\]
by Getzler in \cite{Get1}. Given $m\in\CA$ with $\tilde{m}\circ\tilde{m}=0$, 
he defines
\be\label{thirtyeight}
\tilde{M}_k(\tx_1,\dots,\tx_k)=\left\{
\begin{array}{ll}
0, & k=0 \\
{[\tilde{m},\tx_1]}=\de(\tx_1), & k=1\\
\{ \tilde{m}\}\{ \tx_1,\tx_2\} , & k=2\\
\vdots & \vdots \\
\{ \tilde{m}\}\{ \tx_1,\dots,\tx_k\} & k>1 \\
\vdots & \vdots 
\end{array}\right.\ee
(in a different notation), and proceeds to prove that
\be \tilde{M}\circ\tilde{M}=0\label{thirtynine}\ee
in $\CB$. In short, an $\Ai$ structure $m$ on a graded vector space $A$
is automatically transferred via $M$ to its Hochschild complex $\CA$.
Note that~(\ref{thirtyeight}) generalizes Gerstenhaber's 
construction with $m=m_2$,
and would work equally well with a differential graded associative algebra
where $m=m_1+m_2$. 

\begin{Rm} Recall that the BV bracket is defined by $\pm [m,\tri ]$: higher
BV brackets can then be defined by the above recipe if we replace $\tri$ by
several operators. Any examples from physics already lurking around?
\end{Rm}

\subsection{Strongly homotopy Lie algebras}

A {\bf strongly homotopy Lie algebra} ($\Li$ algebra) is a graded vector
space $A$ plus $n$-ary brackets 
\[ \lbr\; ,\dots,\;\rbr :A^{\ot n}\ra A \]
(one for each $n\geq 1$) satisfying generalizations of the Jacobi identity.
It has been proven in \cite{LM}
that an $\Li$ algebra structure can be obtained from an $\Ai$
structure by super antisymmetrizing $m_n(a_1,\dots,a_n)$'s. In our notation
of coupled braces, we simply set
\be l_n(a_1,\dots,a_n)=
\lbr a_1,\dots,a_n\rbr =\{ m_n\}\{ a_1\}\cdots\{ a_n\},\label{forty}\ee
and observe that the brackets satisfy $|\;\; |$-graded antisymmetry 
(\ref{ds}) as well 
as the higher Jacobi identities 
\be\label{sixfour}\sum_{i+j=n+1}
\sum_{\sigma}(-1)^{p(\sigma ;a_1,\dots,a_n)+i}\;
l_i(l_j(a_{\sigma(1)},\dots,
a_{\sigma(j)}),a_{\sigma(j+1)},\dots,a_{\sigma(n)})=0\;\;\;\forall n,\ee
where $\sigma$ runs through all permutations satisfying
\[ \sigma(1)<\cdots <\sigma(j) \;\;\;\;\mbox{and}\;\;\;\;\sigma(j+1)<\cdots
<\sigma(n).\]

\begin{Rm} There is more than one consistent choice of signs:
the original $\Li$ identities in \cite{LS,LM} involved the
factor $(-1)^{i(j-1)}$, or equivalently $(-1)^{j(i-1)}$, instead of
$(-1)^i$, but the brackets~(\ref{forty}) must satisfy the identities 
in~(\ref{sixfour}) to be consistent with our modified $\Ai$ identities. 
For example, when $n=3$, (\ref{sixfour}) gives
\[ -l_1(l_3(a,b,c))+l_2(l_2(a,b),c)+\cdots,\]
as opposed to
\[ l_1(l_3(a,b,c))+l_2(l_2(a,b),c)+\cdots,\]
which leads to
\[ -m_1(m_3(a,b,c))+m_2(m_2(a,b),c)+\cdots\]
when we replace the $l_i$'s by sums over $m_i$'s. This is the correct 
progression of signs in the $\Ai$ identity~(\ref{thirtyfive}).
\end{Rm}

These identities were not originally conceived as coming from
a master identity like 
\be\label{ai}\{\tilde{m}\circ\tilde{m}\}
\{ sa_1,\dots,sa_n\} =0\;\;\;\forall n\ee 
for $\Ai$ algebras, but from
\[ {\cal C}^2=0\]
where ${\cal C}$ is a coderivation \cite{LS}
(in fact, so did the $\Ai$ identities).
The contribution of the coupled braces notation is to sum up the $\Li$
algebra identities as
\be \label{li}\{\tilde{m}\circ\tilde{m}\}\{ sa_1\}\{ sa_2\}\cdots\{ sa_n\} =0
\;\;\;\forall n.\ee
Writing out the case $m=m_2$, we obtain the super Jacobi identity:
\[ l_2(l_2(a,b),c)+(-1)^{|a|(|b|+|c|)}l_2(l_2(b,c),a)
-(-1)^{|b||c|}l_2(l_2(a,c),b)=0.\]
More generally, we have

\begin{thm}\label{equi}
Identities (\ref{sixfour}) and (\ref{li}) are equivalent.
\end{thm}

\begin{cor}
The existence of the $\Li$ algebra structure~(\ref{li}) 
follows trivially from that of the $\Ai$ algebra structure~(\ref{ai}) by 
taking a sum over all permutations of the $sa_i$.
\end{cor} 

In order to prove the above Theorem, we introduce $\tilde{l}_n$ analogous
to $\tilde{m}_n$: let
\be \tilde{l}_n(a_1,\dots,a_n)\;\defi\; (-1)^{\po}\; l_n(a_1,\dots,a_n).\ee 

\begin{lemma} For every $n\geq 1$, $\tilde{l}_n$ is $\|\;\;\|$-graded
symmetric, that is,
\be \tilde{l}_n(a_{\sigma(1)},\dots,a_{\sigma(n)})={\pttt}\;
\tilde{l}_n(a_1,\dots,a_n).\label{sime}\ee
\end{lemma}

{\it Proof.} See Lemma~\ref{mti}.$\Box$

\begin{prop}\label{pro8} As expected, we have
\be \tilde{l}_n(a_1,\dots,a_n)=\{\tilde{m}_n\}\{ sa_1\}\cdots\{ sa_n\}.\ee
\end{prop}

{\it Proof.} Appendix.$\Box$
\vsp

We will use the two Lemmas below and Proposition~\ref{pi} in the proof of 
Theorem~\ref{equi}.

\begin{lemma} \label{lem1} Given $i+j=n+1$, and a subset $J$ of 
$\{ 1,2,\dots,n\}$ with $|J|=j$, let
\[ \beta_1<\cdots <\beta_j\]
denote the elements of $J$, and
\[ \gamma_1<\cdots <\gamma_{i-1}\]
denote the elements of $J'=\{ 1,2,\dots,n\} \backslash J$. Moreover, fix  
$0\leq k\leq i-1$, and let $\de$, $\dep$ be permutations of the sets $J$
and $J'$ respectively. Finally, if $\sigma_3$ denotes the permutation
\[ \sigma_3=\left( \begin{array}{cccccc}
1 & \dots & j & j+1 & \dots & n \\
\beta_1 & \dots & \bj & \gamma_1 & \dots & \gi \end{array}\right) ,\]
$\sigma_2$ denotes
\[ \sigma_2=\left( \begin{array}{cccccc}
\beta_1 & \dots & \bj & \gamma_1 & \dots & \gi \\
\dbb & \dots & \dbj & \dgb & \dots & \dgi \end{array}\right) ,\]
and $\sigma_1$ denotes
\[\sigma_1=\left( \begin{array}{cccccccccc}
\dbb & \dots&\dots &\dots& \dbj &\dgb & \dots &\dots& \dots&\dgi\\
\dgb & \dots & \dgk &\dbb &\dots& \dots & \dbj & \dgkb &
\dots & \dgi \end{array}\right) ,\]
and if $\si =\sigma_1\sigma_2\sigma_3$, we have
\bea &&\tes\nn\\ &=& \tessu +\tesd +\tesdp +\susb\;\;\;\mbox{(mod 2)} .\nn\eea
\end{lemma}

{\it Proof.} Appendix.$\Box$

\begin{lemma} \label{lem2} With the same notation as above, given $\dpjp$
and $0\leq k\leq i-1$, together with elements $a_0$, $a_{\gamma_1}$, ..., 
$a_{\gi}$
of $A$, let $\sigma_4$ denote the following permutation of $J'\cup \{ 0\}$:
\[ \sigma_4=\left( \begin{array}{cccccccc}
0 & \gamma_1 & \dots & \gamma_{k-1} & \gak & \gakb & \dots & \gki\\
\dgb & \de'(\gamma_2) & \dots & \dgk & 0 & \dgkb & \dots & \dgi
\end{array}\right) .\]
Then we have
\[ \tessd =\tesdp +\susi .\]
\end{lemma}

{\it Proof.} Think of $\sigma_4$ as the identity on 0 and $\dep$ on the rest,
followed by a second permutation which moves 0 to the $(k+1)$-st place.
$\Box$.
\vsp

\begin{prop} \label{pi}Let the notation be as in Lemma~\ref{lem1}. Then 
\be \mms =\sum_{i+j=n+1}\sum_{|J|=j}{\tej}\;\lilj,\ee
where $\sj =\sigma_3$ is 
the permutation sending $a_1$, ..., $a_n$ to $a_{\beta_1}$, ...,
$a_{\bj}$, $a_{\gamma_1}$, ..., $a_{\gi}$.
\end{prop}

{\it Proof.} Appendix.$\Box$

\begin{Rm} It is also possible to write the right hand side as
\bea &&\sum_{i+j=n+1}\sum_J\tilde{\epsilon}(\sigma_5;a_1,\dots,a_n)
(-1)^{\sum_{t=1}^{i-1}\| a_{\gamma_t}\|}\;
\tilde{l}_i(a_{\gamma_1},\dots,a_{\gamma_{i-1}},
\tilde{l}_j(a_{\beta_1},\dots,a_{\beta_j}))\nn\\
&=&\sum_{i+j=n+1}\sum_J\tilde{\epsilon}(\sigma_5;a_1,\dots,a_n)
(-1)^{\sum_{t=1}^{i-1}\| a_{\gamma_t}\|}\;
\{\adj(\tilde{l}_i)\{ a_{\gamma_1},\dots,
a_{\gamma_{i-1}}\}\circ\tilde{l}_j\}\{ a_{\beta_1},\dots,a_{\beta_j}\} ,\nn\eea
where $\sigma_5$ sends $a_1$, ..., $a_n$ to $a_{\gamma_1}$, ..., 
$a_{\gamma_{i-1}}$, $a_{\beta_1}$, ..., $a_{\beta_j}$.
\end{Rm}

The following Proposition settles a natural question:

\begin{prop}\label{bor} We have
\bea &&\lls\nn\\ 
&=&\sum_{i+j=n+1}\sum_{|J|=j}j!(n-j)!(-1)^{\tilde{e}(\sigma_3;a_1,\dots,a_n)}
\;\LiLj ,\nn\eea
hence
\[ \lls =0\]
and
\[ \mms =0\]
are not equivalent! 
\end{prop}

{\it Proof.} Appendix.$\Box$
\vsp

Finally, we have the proof of  equivalence:
\vsp

{\it Proof of Theorem~\ref{equi}.} 
\bea &&\mms\nn\\
&=&\sum_{i,j,J}(-1)^{\tessu}\;\LiLj\nn\\
&=&\sum_{i,j,J}(-1)^{p(\sigma_J;a_1,\dots,a_n)+\sumj +\sumjbt +(i-1)
(1+\sumbt)+\sumig}\nn\\
&&\;\;\;\;\;\;\;\;\;\;\LiLjt\nn\\
&=&\sum_{i,j,J}(-1)^{p(\sigma_J;a_1,\dots,a_n)+\sumj +\sum_{t=1}^j(i+j-1
-t)\| a_{\sigma_J(t)}\| +\sum_{t=j+1}^n(n-t)\| a_{\sigma_J(t)}\| +i-1}\nn\\
&&\;\;\;\;\;\;\;\;\;\;\LiLjt\nn\\
&=&\sum_{i,j,J}(-1)^{p(\sigma_J;a_1,\dots,a_n)+i+\alpha}\;\LiLjt,\nn\eea
where
\[ \alpha =\sum_{t=1}^n(n-t)\| a_t\| -1 \]
is a constant. Note that $(i-1)\|\tilde{l}_j\|$ is $(i-1)$ and not
$(i-1)j$ modulo 2. $\Box$
\vsp

\begin{ex} An $\Li$ algebra structure 
can be imposed on an associative algebra by antisymmetrizing the $\Ai$ 
products described in the previous section.
\end{ex}

\begin{ex} The Hochschild complex $\CB$ of $B=\CA$ is an $\Li$ algebra 
again by virtue of antisymmetrization. 
\end{ex}

\begin{ex} The ``higher-order simple Lie algebras'' introduced by de
Azc\'{a}rraga and Bueno in \cite{AB} are $\Li$~algebras with only one 
higher bracket.
\end{ex}

\begin{ex} See Gnedbaye \cite{Gne} and Hanlon and Wachs \cite{HW}.
\end{ex}

\section{Conclusion}

Coupled braces provide a substantial simplification of the multilinear algebra
of mathematical physics while serving as a stimulant: by analyzing a messy
algebraic relation in terms of compositions, Gerstenhaber brackets, and
higher order differential operators, one usually sees a tidier way of
writing the relation, not to mention several possible ways of generalization
(and conversely, we may be able to compare terms of compact expressions in 
multibraces via simple computer programs and to generate proofs of equalities).
The combination of extensions of maps to derivations and coderivations of
the tensor algebra in one formalism is also fortunate. We hope to continue
exploring this language to unify even more concepts, such as different
types of cohomology theories. Most importantly, we will study the intertwined
homotopy structures on a topological vertex operator algebra -following
Kimura, Voronov, and Zuckerman in~\cite{KVZ}- in an attempt to elevate
them from the shadow of a topological operad to living, breathing algebraic
entities. A master identity for homotopy Gerstenhaber algebras will be
obtained by defining a refinement of multilinear maps and their
compositions \cite{A2}. 
Next, a clear algebraic construction of the $\Ai$ structure
as well as the remaining multibrackets on a TVOA mentioned in~\cite{KVZ} 
should be a top priority (in
addition, the defining identities of general vertex operator algebras 
look familiar in the present context). We expect to
make use of \cite{BDA} which describes relations between the $\Phi$ operators
and $\Li$~algebras.
\vsp

{\it Acknowledgments.} 
Jim Stasheff has been a continual source of support, comments, questions,
references, and people; he is responsible for many little seeds that grew
into viable projects. I would like to thank him especially for mentioning
Kristen Haring's and Lars Kjeseth's theses to me and for the numerous 
technical comments.
I thank Sasha Voronov for the observation that $Hom(TA;TA)$ had been used
by Gerstenhaber and Schack in defining bialgebra cohomology. He has also
been patiently supplying me with information about his joint work with
Kimura and Zuckerman (viva la e-mail!). More thanks to Martin Markl who has
sent me his paper and explained the sign differences.
\pagebreak

\section{Appendix}

{\it Proof of Lemma \ref{lu}.} It suffices to give a proof of this 
statement for a 
transposition $\sigma =\tau =(i,j)$ with $i<j$. We have
\bea &&\sut +\pot \nn\\
&=& \sum_{k=1}^{j-i}(\| a_i\| +\| a_{i+k}\| )+\sum_{k=1}^{j-i-1}(\| a_j\|
+\| a_{i+k}\| )+\pot \nn\\
&=& (j-i)\| a_i\| +\sum_{k=1}^{j-i-1}\| a_{i+k}\| +\| a_j\| +(j-i-1)\| a_j\|
+\sum_{k=1}^{j-i-1}\| a_{i+k}\| +\pot\nn\\
&=& (j-i)\| a_i\| +(j-i)\| a_j\| +\sum_{t\neq i,j}(n-t)\| a_t\| +(n-i)\|
a_j\| +(n-j)\| a_i\|\;\;\;\mbox{(mod 2)}\nn\\
&=& \sum_{t\neq i,j}(n-t)\| a_t\| +(n+j-2i)\| a_j\| +(n-i)\| a_i\|
\;\;\;\mbox{(mod 2)}\nn\\
&=& \po\;\;\;\mbox{(mod 2)}\nn\eea 
(see also Lemma~5 in \cite{Pen1}).$\Box$
\vsp

{\it Proof of Lemma \ref{supsym}.} If $m$ is super antisymmetric, i.e.
\[ \{ m\}\{ a_{\si(1)},\dots,a_{\si(n)}\} =(-1)^{p(\si ;a_1,\dots,a_n)}
\{ m\}\{ a_1,\dots,a_n\} ,\]
then $\tilde{m}$ is $\|\;\|$-graded symmetric, i.e.
\bea &&\{\tilde{m}\}\{ sa_{\si(1)},\dots,sa_{\si(n)}\} \nn\\
&\defi& (-1)^{\sum_{t=1}^n(n-t)\| a_{\si(t)}\|}\{ s\}\{ m\}\{ a_{\si(1)},
\dots,a_{\si(n)}\}\nn\\
&=& (-1)^{\sum_{t=1}^n(n-t)\| a_{\si(t)}\| +p(\si ;a_1,\dots,a_n)}\{ s\}\{ m\}
\{ a_1,\dots,a_n\}\nn\\
&=& (-1)^{\sum_{t=1}^n(n-t)\| a_{\si(t)}\| +\sum_{u<v,\si^{-1}(u)>\si^{-1}
(v)}(\| a_u\| +\| a_v\|)+\tilde{e}(\si ;a_1,\dots,a_n)}\{ s\}\{ m\}\{ a_1,\dots,
a_n\}\nn\\
&=& (-1)^{\sum_{t=1}^n(n-t)\| a_t\| +\tilde{e}(\si ;a_1,\dots,a_n)}\{
s\}\{ m\}\{ a_1,\dots,a_n\}\;\;\;\mbox{by Lemma \ref{lu}}\nn\\
&\defi& (-1)^{\tilde{e}(\si ;a_1,\dots,a_n)}\{\tilde{m}\}\{ sa_1,\dots,
sa_n\} .\Box\nn\eea
\vsp

{\it Proof of Lemma \ref{comut}.} For simplicity let us prove 
this statement for $D(m)=D(x)=2$ and $n=3$. We have
\bea&& {[m,x]}\{ a,b,c\}\nn\\ &=&\{ m\}\{ x\}\{ a,b,c\} -(-1)^{|m||x|+1}\{ x\}
\{ m\}\{ a,b,c\}\nn\\
&=&m(x(a,b),c)+(-1)^{|a||x|-1}m(a,x(b,c))+(-1)^{|m||x|}(x(m(a,b),c)+
(-1)^{|m||a|-1}x(a,m(b,c)))\nn\\
&=&m(x(a,b),c)-(-1)^{|a||x|}m(a,x(b,c))+(-1)^{|m||x|}x(m(a,b),c)
-(-1)^{|m||a|+|m||x|}x(a,m(b,c)).\nn\eea
On the other hand,
\bea &&{[\tm,\tx ]}\{ sa,sb,sc\}\nn\\
&=&\{\tm\}\{\tx\}\{ sa,sb,sc\} -(-1)^{\| m\|\;\| x\|}\{\tx\}\{\tm\}\{
sa,sb,sc\}\nn\\
&=&\tm(\tx(a,b),c)+(-1)^{\| a\|\;\| x\|}\tm(a,\tx(b,c))-(-1)^{\| m\|\;\| x\|}
(\tx(\tm(a,b),c)+(-1)^{\| m\|\;\| a\|}\tx(a,\tm(b,c)))\nn\\
&=&(-1)^{\| a\| +\| x\| +\| a\| +\| b\|}m(x(a,b),c)+(-1)^{\| a\|\;\| x\|
+\| a\| +\| b\|}m(a,x(b,c))\nn\\
&&+(-1)^{\| m\|\;\| x\| +\| a\| +\| m\| +\| a\| +\| b\| +1}x(m(a,b),c)
+(-1)^{\| m\|\;\| a\| +\| m\|\;\| x\| +1+\| a\| +\| b\|}x(a,m(b,c))\nn\\
&=&(-1)^{\| x\| +\| b\|}(m(x(a,b),c)-(-1)^{(\| a\| +1)(\| x\| +1)}
m(a,x(b,c))\nn\\
&&+(-1)^{(\| m\| +1)(\| x\| +1)}x(m(a,b),c)-(-1)^{(|m|+1)
(\| a\| +\| x\|)+\| a\| +\| x\|}x(a,m(b,c))),\nn\eea
which is the right hand side.$\Box$
\vsp

{\it Proof of Lemma \ref{B}.} 
\bea &&{[x,y\cdot z]}-[x,y]\cdot z-(-1)^{d(x)(d(y)+1)}y\cdot [x,z]\nn\\
&=&\{ x\}\{ y,z\} -(-1)^{d(x)d(y\cdot z)}\{ y\cdot z\}\{ x\}
-\{ x\}\{ y\}\cdot z\nn\\&&+(-1)^{d(x)d(y)}\{ y\}\{ x\}\cdot z
-(-1)^{d(x)d(y)+d(x)}y\cdot\{ x\}\{ z\} +(-1)^{d(x)d(y)+d(x)+d(x)d(z)}y
\cdot\{ z\}\{ x\}\nn\\
&=&\{ x\}\{ y\cdot z\}\nn\\
&&-(-1)^{d(x)(d(y)+d(z)+1)}(y\cdot\{ z\}\{ x\} +(-1)^{d(x)(d(z)+1)}
\{ y\}\{ x\}\cdot z)\nn\\
&&-\{ x\}\{ y\}\cdot z+(-1)^{d(x)d(y)}\{ y\}\{ x\}\cdot z\nn\\
&&-(-1)^{d(x)d(y)+d(x)}y\cdot\{ x\}\{ z\}+(-1)^{d(x)(d(y)+d(z)+1)}y
\cdot\{ z\}\{ x\}\nn\eea
by (\ref{twentyseven}) where $n=1$, $x_1=y$, $x_2=z$, and $y_1=x$. After
cancellations, we obtain
\bea &&\{ x\}\{ y\cdot z\} -\{ x\}\{ y\}\cdot z-(-1)^{d(x)D(y)}y\cdot
\{ x\}\{ z\}\nn\\
&=&(-1)^{d(x)+D(y)}(\de(\{ x\}\{ y,z\})-\{\de(x)\}\{ y,z\}
-(-1)^{d(x)}\{ x\}\{\de(y),z\} -(-1)^{d(x)+d(y)}\{ x\}\{
y,\de(z)\})\nn\eea
by (\ref{twentyeight}) with $n=1$, $y_1=y$, and $y_2=z$.$\Box$
\vsp

{\it Proof of Proposition \ref{eqq}.} First of all, we have
\[ {[m_i,m_j]}=m_i\circ m_j-(-1)^{(i-1)(j-1)+ij}m_j\circ m_i=m_i\circ m_j
-(-1)^nm_j\circ m_i.\]
Then
\bea &&\sum_{i+j=n+1}(-1)^{j}[m_i,m_j]\{ a_1,\dots,a_n\}\nn\\
&=&\sum_{i+j=n+1}(-1)^j\{ m_i\circ m_j\}\{ a_1,\dots,a_n\} +\sum_{i+j=n+1}
(-1)^{j+n+1}\{ m_j\circ m_i\}\{ a_1,\dots,a_n\}\nn\\
&=&\sum_{i+j=n+1}(-1)^j\{ m_i\circ m_j\}\{ a_1,\dots,a_n\} +\sum_{i+j=n+1}
(-1)^i\{ m_j\circ m_i\}\{ a_1,\dots,a_n\}\nn\\
&=&2\sum_{i+j=n+1}(-1)^j\{ m_i\circ m_j\}\{ a_1,\dots,a_n\}\nn\\
&=&2\sum_{i+j=n+1}(-1)^j\sum_{k=0}^{i-1}(-1)^{j(k-1)+j(|a_1|+\cdots +
|a_k|)}m_i(a_1,\dots,a_k,m_j(a_{k+1},\dots,a_{k+j}),\dots,a_n)\nn\\
&=&2\sum_{i+j=n+1}\sum_{k=0}^{i-1}(-1)^{j(|a_1|+\cdots +|a_k|)+jk+j+k}
m_i(a_1,\dots,a_k,m_j(a_{k+1},\dots,a_{k+j}),\dots,a_n). \Box\nn\eea
\vsp

{\it Proof of Proposition \ref{pro8}.} We have
\bea &&\{\tilde{m}_n\}\{ sa_1\}\cdots\{ sa_n\}\nn\\&=&\sum_{\sigma}{\pttt}\;
\tilde{m}_n(a_{\sigma(1)},\dots,a_{\sigma(n)})\nn\\
&=&\sum_{\sigma}(-1)^{\p +\su +\pos}\; 
m_n(a_{\sigma(1)},\dots,a_{\sigma(n)})\nn\\
&=&\sum_{\sigma}(-1)^{\p +\po}\; m_n(a_{\sigma(1)},\dots,a_{\sigma(n)})
\;\;\;\mbox{by Lemma~\ref{lu}}\nn\\
&=&(-1)^{\po}\; l_n(a_1,\dots,a_n)\nn\\
&=& \tilde{l}_n(a_1,\dots,a_n).\Box\nn\eea 
\vsp

{\it Proof of Lemma \ref{lem1}.} By repeated 
applications of (\ref{split}), we get (modulo 2)
\bea &&\tebiu \nn\\ &=& \tessu +\tessi +\tessb \nn\\
&=& \tessu +\tesd +\tesdp +\susa \nn\\&=&\tessu +\tesd +\tesdp
+\susb.\Box\nn\eea
\vsp

{\it Proof of Proposition \ref{pi}.} Let $\si$ be as in Lemma~\ref{lem1}. Then
\bea &&\mms\nn\\ &=& \sum_{i+j=n+1}\sum_{k=0}^{i-1}\sum_J\sum_{\dpj}
\sum_{\dpjp}(-1)^{\tes +\susu}\nn\\ &&\;\;\;\;\;\;\;\;\;\;\mimj\nn\\
&=& \sum_{i,j,k,J,\de,\dep}(-1)^{\tessu +\tesd 
+\tesdp +\susd}\nn\\ &&\;\;\;\;\;\;\;\;\;\;\mimjt\nn\\ 
&=& \sum_{i,j,k,J,\dep}(-1)^{\tessu +\tesdp +\susd}\nn\\ 
&&\;\;\;\;\;\;\;\;\;\;\milj .\nn\eea
We would also like to eliminate the sums over $k$ and $\dep$ and change 
$\tilde{m}_i$ into $\tilde{l}_i$ by introducing the appropriate 
symmetrization factor. By Lemma~\ref{lem2}, this factor is exactly
\[ (-1)^{\tesdp +\susd}.\]
Then we have
\be \mms =\sum_{i,j}\sum_J(-1)^{\tessu}\;\lilj ,\ee
where $\sigma_3 =\sj$ is the permutation that sends $a_1$, ..., $a_n$ into the
sequence on the right hand side. $\Box$
\vsp

{\it Proof of Proposition \ref{bor}.} We follow the notation and 
methods of the previous proof.
Proceeding in exactly the same manner, we obtain
\bea &&\lls\nn\\
&=&\sum_{i,j,k,J,\de,\dep}(-1)^{\tessu +\tesd +\tesdp +\susd}\nn\\
&&\;\;\;\;\;\;\;\;\;\;\liljt\nn\\
&=&\sum_{i,j,k,J,\dep}j!(-1)^{\tessu +\tesdp +\susd}\nn\\
&&\;\;\;\;\;\;\;\;\;\;\lilju\nn\\
&=&\sum_{i,j,J}j!(n-j)!(-1)^{\tessu}\;\LiLj\nn\eea
as each $\tilde{l}_i$ is suspended-graded symmetric. $\Box$
\vsp


\begin{thebibliography}{99}

\bibitem{A} F. AKMAN, On some generalizations of Batalin-Vilkovisky
algebras, to appear in {\em JPAA}, preprint q-alg/9506027.

\bibitem{A2} F. AKMAN, A master identity for homotopy Gerstenhaber
algebras, in preparation.

\bibitem{AB} J.A. DE AZC\'{A}RRAGA AND J.C.P. BUENO, Higher-order simple Lie
algebras, {\it Commun. Math. Phys.} {\bf 184} (1997), 669-681;
preprint hep-th/9605213.

\bibitem{BV} I.A. BATALIN AND G.A. VILKOVISKY, Gauge algebra and quantization,
{\em Phys. Lett.} {\bf 102B} (1981), 27-31. 

\bibitem{BDA} K. BERING, P.H. DAMGAARD, AND J. ALFARO, Algebra of higher
antibrackets, {\it Nucl. Phys.} {\bf B478} (1996), 459,
preprint hep-th/9604027.

\bibitem{FN} A. FR\"{O}LICHER AND A. NIJENHUIS, Some new cohomology 
invariants for complex manifolds (I and II), {\it Koninklijke Nederlandse
Akademie van Wetenshappen (later Indag. Math.), Series A, Proceedings} 
{\bf 59} (1956), 540.

\bibitem{Ger} M. GERSTENHABER, The cohomology structure of an associative ring,
{\em Ann. Math.} {\bf 78} (1963), 267-288.

\bibitem{GeS} M. GERSTENHABER AND S.D. SCHACK, Bialgebra cohomology,
deformations, and quantum groups, {\it Proc. Nat. Acad. Sci.} {\bf 87}
(1990), 478-481.

\bibitem{Get1} E. GETZLER, Cartan homotopy formulas and the Gauss-Manin
connection in cyclic homology, in "Quantum deformations of algebras and
their representations", {\em Israel Mathematical Conference
Proceedings}, Vol. 7, 1993 (eds. A. Joseph and S. Shnider).

\bibitem{GS} E. GETZLER AND J.D.S. JONES, $A_{\infty}$-algebras and the
cyclic bar complex, {\em Illinois J. Math.} {\bf 34} (1989), 256-283.

\bibitem{Gia} A. GIAQUINTO, Deformation methods in quantum groups, Ph.D.
Thesis, University of Pennsylvania, 1991.

\bibitem{Gne} V. GNEDBAYE, Operads of k-ary algebras, in "Operads: Proceedings
of Renaissance Conferences", J.-L. Loday, J.~Stasheff, and A.A. Voronov,
eds., {\em Contemp. Math.}, vol. 202, AMS, Providence, RI (1996).

\bibitem{HW} P. HANLON AND M.L. WACHS, On Lie-k algebras, {\em Adv. in Math.}
{\bf 113} (1995), 206-236.

\bibitem{Ha} K.A. HARING, On the events leading to the formulation of the 
Gerstenhaber algebra: 1945-1966, Master's Thesis, University of North 
Carolina at Chapel Hill, 1995.

\bibitem{KVZ} T. KIMURA, A.A. VORONOV, AND G.J. ZUCKERMAN, Homotopy
Gerstenhaber algebras and topological field theory, in
"Operads: Proceedings of Renaissance Conferences", J.-L. Loday,
J.~Stasheff, and A.A. Voronov, eds., {\em Contemp. Math.}, vol. 202,
AMS, Providence, RI (1996); preprint q-alg/9602009.

\bibitem{Kje} L. KJESETH, BRST cohomology and homotopy Lie-Rinehart pairs,
Ph.D. Thesis, University of North Carolina at Chapel Hill, 1996.

\bibitem{Ko} J.-L. KOSZUL, Crochet de Schouten-Nijenhuis et cohomologie,
{\em Ast\'{e}risque} (1985), 257-271.

\bibitem{LM} T. LADA AND M. MARKL, Strongly homotopy Lie algebras, {\it
Communications in Algebra} {\bf 23} (1995), 2147-2161; preprint
hep-th/9406095.

\bibitem{LS} T. LADA AND J.D. STASHEFF, Introduction to sh Lie algebras for
physicists, {\it Intern'l J. Theor. Phys.} {\bf 32} (1993), 1087-1103,
preprint hep-th/9209099, UNC-MATH-92/2.

\bibitem{LZ} B.H. LIAN AND G.J. ZUCKERMAN, New perspectives on the
BRST-algebraic structure of string theory, {\em Commun. Math. Phys.} 
{\bf 154} (1993), 613-646; preprint hep-th/9211072; MR 94e:81333.

\bibitem{Ma} M. MARKL, A cohomology theory for $A(m)$-algebras and
applications, {\it JPAA} {\bf 83} (1992), 141-175.

\bibitem{Pen1} M. PENKAVA, $L_{\infty}$ algebras and their cohomology,
preprint q-alg/9512014.

\bibitem{Pen2} M. PENKAVA, Infinity algebras and the homology of graph
complexes, preprint q-alg/9601018.

\bibitem{PS} M. PENKAVA AND A. SCHWARZ, $A_{\infty}$ algebras and the
cohomology of moduli spaces, preprint hep-th/9408064.

\bibitem{St} J.D. STASHEFF, Homotopy associativity of {\em H}-spaces I and
II, {\em AMS Trans.} {\bf 108} (1963), 275-292 and 293-312.

\bibitem{Jim} J.D. STASHEFF, The intrinsic bracket on the deformation
complex of an associative algebra, to appear in {\it JPAA}; preprint
UNC-MATH-91/1.

\bibitem{VG} A.A. VORONOV AND M. GERSTENHABER, Higher order operations
on Hochschild complex, {\em Functional Anal. Appl.} {\bf 29} (1995), no.1,
1-6.

\end{thebibliography}
\end{document}